\documentstyle[12pt]{article}
\def\m@th{\mathsurround=0pt }
\def\eqalign#1{\null\,\vcenter{\openup\jot \m@th
   \ialign{\strut\hfil$\displaystyle{##}$&$
      \displaystyle{{}##}$\hfil \crcr#1\crcr}}\,}
\newtheorem{theorem}{Theorem}[section]

\newtheorem{corollary}[theorem]{Corollary}
\newtheorem{lemma}[theorem]{Lemma}

\newtheorem{definition}[theorem]{Definition}

\def\ad{{\rm ad }\ }

\def\adr{{\rm ad}_r\ }

\font\goth=eufm10

\def\gg{\hbox{\goth g}}

\def\gn{\hbox{\goth n}}

\def\gh{\hbox{\goth h}}

\def\ga{\hbox{\goth a}}

\def\gr{\hbox{\goth r}}

\def\tip{\hbox{\rm tip}}

\begin{document}

\title
{  Invariant Differential Operators for Quantum Symmetric 
Spaces, I  }
\author { Gail Letzter\thanks{supported by NSA grant no. 
MDA904-03-1-0033.
AMS subject classification 17B37}\\
Mathematics Department\\
Virginia Tech\\
Blacksburg, VA 24061\\
letzter@math.vt.edu}

\maketitle
 
\begin{abstract} This is the first paper in   a series of two 
which proves a  
version of a theorem of  Harish-Chandra for quantum symmetric spaces in the maximally split case:
 There is a  Harish-Chandra 
map which induces an isomorphism  between the 
  ring of  quantum invariant differential 
operators   and the ring of invariants of 
a certain   
Laurent polynomial ring  under an action of the 
restricted Weyl group.  Here, we establish this result for 
all quantum symmetric spaces defined using 
irreducible symmetric pairs not of type  EIII, EIV, EVII, or EIX.  A quantum version of a 
related   theorem  due to Helgason is also obtained: The image of the 
center under this Harish-Chandra map is the entire invariant ring   if 
and only if the underlying irreducible 
symmetric pair is not one of these four types. 
    \end{abstract}

\section* {Introduction}

In [HC2, Section 4], Harish-Chandra proved the following fundamental result for 
semisimple Lie groups: the Harish-Chandra map induces
 an isomorphism between
the ring of invariant differential operators on a symmetric space and 
invariants of an appropriate polynomial ring under the restricted Weyl group. When the 
symmetric space is simply a complex semisimple Lie group, this 
result
is Harish-Chandra's famous realization of the center of the enveloping 
algebra of a semisimple Lie algebra as Weyl group invariants of the Cartan 
subalgebra ([HC1]).   Helgason further refined Harish-Chandra's result by 
analyzing
the image  of those invariant differential operators which come from the center of 
the enveloping algebra under the Harish-Chandra 
homomorphism ([He]).  This paper is the first of   two papers 
which   establish quantum analogs for these results of  
Harish-Chandra and Helgason in  the recently developed theory of quantum symmetric spaces
 (see   [L3],[L4], and [L5]).

There is an underlying symmetric pair of Lie algebras 
$\gg,\gg^{\theta}$ associated to each symmetric space. Here, $\gg$ is a complex semisimple 
  Lie algebra and 
  $\gg^{\theta}$ is the Lie subalgebra fixed by an involution 
  $\theta$.     
In this algebraic framework,   invariant differential operators on 
the symmetric space
 correspond to  $\ad \gg^{\theta}$ invariant elements of 
  $U(\gg)$. Thus quantum invariant differential operators should 
  ``be'' elements of the quantum analog of 
 the fixed ring of $U(\gg)$ under the action of $\gg^{\theta}$.  
 In fact, a quantum symmetric pair consists of the quantized enveloping algebra 
of $\gg$ and a    coideal subalgebra $B$ which plays the role of
 $U(\gg^{\theta})$. The study of quantum invariant differential 
  operators is an analysis of the ring $\check U^B$ of invariants 
  in  the (simply connected) quantized enveloping algebra $\check U$  
  of $\gg$  
 with respect to  the right adjoint action of $B$.  (We assume that 
 $\check U$ is an algebra over the algebraic closure ${\cal C}$ of the 
 rational function field ${\bf C}(q)$.)
  
 Let $\ga$ denote the eigenspace for $\theta$ with eigenvalue $-1$ 
  inside a maximally split Cartan subalgebra $\gh$ of $\gg$. 
  In the classical case, the Harish-Chandra map associated to 
 $\gg,\gg^{\theta}$ is the projection from 
  $U(\gg)$ onto the enveloping algebra of $\ga$, defined using the 
 Iwasawa decomposition of $\gg$. 
  This picture can be lifted to the quantum setting as follows.
 The Cartan subalgebra $\check U^0$ of $\check U$ is the group 
  algebra of a multiplicative group $\check T$  isomorphic to the weight lattice 
  of $\gg$. The quantum analog of $\ga$ is the 
  multiplicative subgroup    $\check {\cal A}$ of  an extension of $\check T$  
  isomorphic to the image of the weight lattice of $\gg$ under the restriction 
  map from $\gh^*$ to $\ga^*$.  Using a quantum version of the 
  Iwasawa decomposition (Theorem 2.2 below), the quantum 
  Harish-Chandra map ${\cal P}_B$ is a projection of $\check U$ onto 
    ${\cal C}[\check {\cal A}]$. 
    
    The restricted root system $\Sigma$ associated to the pair 
    $\gg,\gg^{\theta}$ spans the vector space $\ga^*$.  Thus the 
   action of the  restricted Weyl group $W_{\Theta}$ on $\Sigma$ 
   induces an action of $W_{\Theta}$ on $\check {\cal 
   A}$.  This action can be further twisted to a dotted (or 
   translated) action of $W_{\Theta}$ on ${\cal C}[\check {\cal A}]$.  
   Let ${\cal A}$ denote the subgroup of $\check {\cal A}$ 
   corresponding to the weight lattice of the root system $2\Sigma$. 
    The 
  center of $\check U$, denoted by $Z(\check U)$, is a subring of 
  $\check U^B$. 
  
  Recall that irreducible symmetric pairs 
  $\gg,\gg^{\theta}$ are classified in [A], and each pair or family of pairs 
  is given a name which is called its ``type''. We assume that 
  $\gg,\gg^{\theta}$ is irreducible.   The analogous results for 
  arbitrary symmetric pairs of Lie algebras are easily deduced from 
  their decomposition into a direct sum of irreducible ones. The main result of this paper is

 \medskip
 \noindent
 {\bf Theorem:}{\it (see Theorem 6.1) The image of $Z(\check U)$ under ${\cal P}_B$ is 
 equal to ${\cal C}[{\cal 
 A}]^{W_{\Theta}.}$ if and only 
 if $\gg,\gg^{\theta}$ is not of
  type EIII, EIV, EVII, or EIX.   
 Moreover, in these cases, ${\cal P}_B(\check U^B)={\cal P}_B(Z(\check U)).$ }  
 \medskip

 Perhaps the most intricate part in the proof of the above theorem is showing that 
 ${\cal P}_B(\check U^B)$ is invariant under the dotted action of 
 $W_{\Theta}$.    Recall that there is no direct quantum version of 
  the Lie group $G$ associated to $\gg$; there is only a   quantum analog of the ring of 
  regular functions on $G$.  
   Thus the analytic techniques 
used in Harish-Chandra's original work [HC2] are unavailable here.   
Lepowsky's  algebraic proof [Le] 
using  the Cartan decomposition of $\gg$ with respect to 
$\gg^{\theta}$  cannot be adapted to the 
quantum setting because   there is no obvious quantum Cartan 
decomposition.  For the cases considered in the theorem, the dotted 
$W_{\Theta}$ invariance of ${\cal P}_B(\check U^B)$  ultimately follows 
from the dotted $W_{\Theta}$ invariance of ${\cal P}_B(Z(\check U))$ and the 
equality of these two algebras established at the end of the paper. 

   The first step in the proof of the   theorem  is to refine
the codomain of ${\cal P}_B$ upon restriction to $\check U^B$.   
  We exploit the fact  that $\check U^B$ is a subring of   the locally 
finite part $F_r(\check U)$ of $\check U$ in order to show
  ${\cal P}_B(\check U^B)$ is contained in   ${\cal C}[{\cal 
  A}]$.  
Using a filtration of $\check U$ introduced in [L5, Section 5], we further show that the 
possible highest degree terms of elements in ${\cal P}_B(\check U^B)$ 
are linear combinations of elements in ${\cal A}$ which correspond to antidominant integral weights in the weight lattice of 
$2\Sigma$.

It turns out that the restriction 
of ${\cal P}_B$ to $Z(\check U)$ is just the 
composition of the ordinary quantum Harish-Chandra map ${\cal P}$ 
followed by restriction of ${\cal P}_B$ to  $\check U^0$.
It is well known that ${\cal P}(Z(\check U))$ consists of the 
invariants of a particular Laurent polynomial subring of $\check U^0$ with respect to the dotted action of the large Weyl 
group. Thus determining ${\cal P}_B(Z(\check U))$ 
reduces to computing the image  of the invariants of this Laurent 
polynomial ring inside ${\cal 
C}[{\cal A}]$.  
When $\gg,\gg^{\theta}$ is one of the four types EIII, EIV, EVII, EIX, 
we show that  ${\cal P}_B(Z(\check U))\neq {\cal 
C}[{\cal A}]^{W_{\Theta}.} $ using specialization at $q=1$ and   the 
corresponding classical 
results.   Unfortunately, specialization fails to show the other
direction of the primary assertion of 
the  theorem.  Instead, we use 
information relating the 
 the weight lattice  of  the root 
system of $\gg$ to that of  $\Sigma$ and a special basis 
of $Z(\check U)$ to determine ${\cal P}_B(Z(\check U))$.   The proof 
turns out to be more delicate for those      
 symmetric pairs related to the irreducible symmetric pair of type AII. 
 In Section 6, we use combinatorial and character formula  arguments 
 to handle this special family. 
  The last assertion of the theorem follows from  a comparison of the set of highest degree terms of each of the  
 three algebras    ${\cal P}_B(Z(\check U))$, ${\cal 
 C}[{\cal A}]^{W_{\Theta}.}$, and ${\cal P}_B(\check U^B)$.
 
 Now suppose that $\gg,\gg^{\theta}$ is of type EIII, EIV, EVII, or 
 EIX.  In the sequel to this paper [L6],  we show that ${\cal P}_B(\check U^B)={\cal C}[{\cal 
 A}]^{W_{\Theta}.}$ in these cases as well. Thus  this paper
 and [L6]  establish complete   quantum analogs of   Harish-Chandra's 
 and Helgason's results on the image of invariant differential operators  
 under the Harish-Chandra map.  Moreover, in [L6], we find a 
 particular nice basis for ${\cal P}_B(\check U^B)$ related to 
a special family of   difference operators for Macdonald 
 polynomials.

There is  another function on $\check U$, the radial component map 
${\cal X}$,  which is closely related to ${\cal P}_B$ and is useful in 
  analyzing ${\cal P}_B(\check U^B)$.  We show that the map which 
sends ${\cal P}_B(u)$ to ${\cal X}(u)$  defines an  algebra isomorphism 
of ${\cal P}_B(\check U^B)$ onto ${\cal X}(\check U^B)$. This enables
us to determine the kernel of ${\cal P}_B(\check U^B)$.  It should be 
noted that the quotient 
of $\check U^B$ modulo this kernel is   the exact quantum analog of the 
ring of invariant differential operators on the corresponding
symmetric space. 
Furthermore, the main theorem and [L6, Corollary 4.2] ensures that ${\cal X}(\check U^B)$ is a 
polynomial ring in rank $\Sigma$ variables.  With the help of ${\cal 
X}$, we interpret elements of 
$\check U^B$ as difference operators acting on the character group 
ring of $2\Sigma$.  It is precisely this interpretation that leads 
  to the identification of quantum zonal spherical functions as
Macdonald (or  Macdonald-Koornwinder) polynomials 
  (see for example [N], [NS], [L5],  [DN],   [NDS], and [DS]).    The   theorem and 
its extension in [L6] ensure that 
  ${\cal 
X}(\check U^B)$  forms a completely integrable system  of difference operators 
whose eigenfunctions are the  orthogonal
polynomials associated to quantum zonal spherical 
functions.

We briefly describe the organization of this paper.   Section 1 sets 
notation and recalls basic facts about quantized enveloping algebras, 
restricted root systems, and quantum symmetric pairs.  In Section 2, 
we review and generalize  the construction of the 
Harish-Chandra map ${\cal P}_B$ defined in [L5]. Section 3 
is a detailed study of the image  of the weight lattice 
$P(\pi)$ under the map $\ \tilde{\ }$  defined by  restriction to $\ga$. 
 The 
  definition of the quantum 
  radial component map ${\cal X}$ and the filtration on $\check U$
  as described in [L5]  are   extended to the general case in Section 4.
 Connections between ${\cal P}_B(\check U^B)$ and ${\cal X}(\check 
  U^B)$ are established, leading to   a nice description of the possible highest 
  degree terms of elements in ${\cal P}_B(\check U^B)$. The main
  result   is proved in Sections 5 and 6.    Section 5 handles the 
  four exceptional cases EIII, EIV, EVII, and EIX  using specialization and those irreducible 
  symmetric pairs which satisfy $\widetilde{P^+(\pi)}=P^+(\Sigma)$.
    Section 6 
  is devoted to the proof of the main theorem for those 
   irreducible symmetric pairs $\gg,\gg^{\theta}$ which contain a 
  symmetric pair of type AII.

  \medskip
  \noindent
  {\bf Acknowledgements.}  The author would like to thank Nolan 
  Wallach for his wonderful insight and suggestions which led to this 
  project and Dan Farkas for his helpful comments on writing.

    \section{Background and Notation }
    Let ${\bf C}$ denote the complex numbers, ${\bf Q}$ denote the 
rational numbers, ${\bf Z}$ denote the integers,  ${\bf R}$ denote 
the real numbers, ${\bf N}$ denote the nonnegative integers, and $q$ 
denote an indeterminate. 
Let $\{q^r|\ r\in{\bf Q}\}$ denote the multiplicative group 
isomorphic to ${\bf Q}$ using the map $q^r\mapsto r$.  
Write ${\cal C}$ for the algebraic closure of  the group algebra of $\{q^r|\ r\in{\bf Q}\}$
over ${\bf C}$.  Note that the group algebra of $\{q^r|\ r\in{\bf Q}\}$
over ${\bf R}$ can be made into an ordered field in a manner similar 
to the rational function field ${\bf R}(q)$ (see [Ja, Section 5.1]). Let 
${\cal R}$ be the real algebraic closure of the group algebra of $\{q^r|\ r\in{\bf Q}\}$
over ${\bf R}$.  

Suppose that  $\Phi$ is a root system.    Write $Q(\Phi)$ for the root lattice and let 
$Q^+(\Phi)$ be the subset of $Q(\Phi)$ equal to the ${\bf N}$ span 
of the positive roots in $\Phi$.     Let $P(\Phi)$ denote the weight lattice of $\Phi$ and 
  let $P^+(\Phi)$ be the subset of $P(\Phi)$ consisting of 
dominant integral weights.   (Sometimes we replace $\Phi$ with the 
symbol used to represent the positive simple roots in the notation for 
the weight and root lattices and their subsets.)

Let $\gg$ be a complex semisimple Lie algebra with triangular 
decomposition $\gg=\gn^-\oplus \gh\oplus \gn^+$. 
Let $\theta$ be a maximally split 
involution with respect to the fixed Cartan subalgebra $\gh$ of $\gg$ 
and the above triangular decomposition in the 
sense of [L3, (7.1), (7.2), and (7.3)]. Write $\gg^{\theta}$ for the fixed Lie subalgebra of 
$\gg$ with respect to $\theta$. The pair $\gg,\gg^{\theta}$ is a 
classical (infinitesimal) symmetric pair.  
We assume throughout the paper that $\gg,\gg^{\theta}$ is an irreducible 
symmetric pair  in the sense of [A] (see also [L4, Section 7]). The results of this paper extend in a 
straightforward manner to the   
general case.  

There are two root systems 
associated to $\gg,\gg^{\theta}$.    The first is the ordinary root 
system of $\gg$ which we denote by $\Delta$. Let $(\ , \ )$ denote the 
corresponding Cartan inner product on $\Delta$. Here we assume that 
the positive roots of $\Delta$ 
correspond to the root vectors in $\gn^+$.  Write  
$\pi=\{\alpha_1,\dots, \alpha_n\}$ for the positive simple roots of 
$\Delta$.  We use ``$\leq $'' to denote the usual partial order on $\gh^*$.   In 
particular, given $\alpha$ and $\beta$ in $\gh^*$, we say that 
$\alpha\leq\beta$ if and only if $\beta-\alpha\in Q^+(\pi)$.  

The second root system is the set of restricted roots 
$\Sigma$ defined using the involution $\theta$.   More precisely,  
  $\theta$ induces an involution $\Theta$ on $\gh^*$ which restricts 
  to an involution on  
  $\Delta$. 
Given $\alpha\in \gh^*$, set $\tilde\alpha= (\alpha-\Theta(\alpha))/2$.  
The restricted root system  $\Sigma$ is the set 
$$\Sigma=\{\tilde\alpha|\alpha\in \Delta{\rm \  and \ }\Theta(\alpha)\neq
\alpha\}.$$ Moreover, $\Sigma$   inherits the structure of a root system using  
the inner product of $\Delta$ ([Kn, Chapter VI, Section 4]).  
Recall [L3 ,Section 7, (7.5)] that there exists a permutation ${\rm p}$ of the set 
$\{1,2,\dots, n\}$ such that ${\rm p}$ induces a diagram automorphism on $\pi$ 
and $$\Theta(-\alpha_i)-\alpha_{{\rm p}(i)}\in 
Q^+(\pi_{\Theta})\eqno{(1.1)}$$ for each $i$. Set
$\pi^*=\{\alpha_i\in \pi\setminus \pi^*|\ i\leq {\rm p}(i)\}$. 
The set $\{\tilde\alpha_i|\ \alpha_i\in \pi^*\}$ is the set
of positive simple roots for the root system $\Sigma$.

  Let $U=U_q(\gg)$ denote the quantized enveloping algebra of $\gg$. 
  The algebra $U$ is a Hopf algebra over ${\cal C}$ with generators
  $x_i,y_i,t_i^{\pm 1}$, $1\leq i\leq 
n$,  subject to the relations and Hopf structure given  in [L3, Section 1, (1.4)-(1.10)] or [Jo, 
3.2.9].   Let 
$T$ denote the group generated by $t_i$,  $1\leq i\leq n$, let $U^+$
denote the subalgebra of $U$ generated by $x_i, 1\leq i\leq n$, and let
$G^-$ denote the subalgebra of $U$ generated by $y_it_i, 1\leq i\leq n$.
 Let $U^0$ denote the group algebra generated by $T$.

Let $\tau$ denote the  group isomorphism from the additive group 
$Q(\pi)$ to the multiplicative group $T$ defined by 
$\tau(\alpha_i)=t_i$ for $1\leq i\leq n$. 
Given a vector subspace $S$ of $U$, the $\beta$ weight space of $S$, 
denoted by  $S_{\beta}$, is  the set 
all elements $s$ in $S$ which satisfy 
$$\tau(\gamma)s\tau(-\gamma)=q^{(\gamma,\beta)}s \eqno{(1.2)}$$ 
for all 
$\tau(\gamma)\in T$.

Sometimes it will be necessary to consider larger algebras than     $ U$ 
which are obtained using extensions of  the group $T$.   In particular,  suppose that  $M$ is a multiplicative 
monoid isomorphic to an  additive submonoid of    $\sum_{1\leq 
i\leq n}{\bf Q}\alpha_i$ via the obvious extension of $\tau$.  Then the 
algebra $ UM$ 
is just the algebra generated by $  U$ and $M$ 
subject to the additional relation (1.2) for each $\tau(\gamma)\in M$ 
and $s\in U_{\beta}$.   It 
is straightforward to check that when $M$ is a group,   the Hopf algebra structure of $U$ 
extends to 
$ UM$. Recall that the augmentation ideal of a Hopf algebra is just 
the kernel of the counit. If $UM$ is a Hopf algebgra and $A$ is a 
subalgebra, then we write $A_+$ to denote the intersection of a 
  $A$  with the augmentation ideal of $AM$.

The most common  extension of the type described in the previous paragraph is the simply connecting 
quantized enveloping algebra, denoted  by $\check U$ ([Jo, 
Section 3.2.10]).  Here, the multiplicative monoid $M$ is the group 
$\check T$ which is the extension of $T$ isomorphic to the weight 
lattice $P(\pi)$ via $\tau$. 
Let
$\check U^0$ denote the group algebra generated by $\check T$.

  A   quantum analog  of $U(\gg^{\theta})$ 
inside of $U_q(\gg)$ is  a  maximal left coideal 
subalgebra  of $U_q(\gg)$ which specializes to $U(\gg^{\theta})$ as $q$ 
goes to $1$.   Quantum analogs of $U(\gg^{\theta})$ 
inside of $U_q(\gg)$ are classified in [L3, Section 7]; a full 
description of the  generators and relations for all quantum analogs 
of $U(\gg^{\theta})$ 
associated to each irreducible symmetric pair $\gg,\gg^{\theta}$ can be 
found in [L4, Section 7].
We briefly describe here the construction of these  coideal subalgebras.     Set
$\pi_{\Theta}=\{\alpha_i|\Theta(\alpha_i)=\alpha_i\}$.  Let ${\cal M}$ denote the Hopf 
subalgebra  of $U$ generated by $x_i,y_i, t_i^{\pm 1}$ for 
$\alpha_i\in \pi_{\Theta}$.       Let $T_{\Theta}$ be the subgroup of 
$T$ defined by $T_{\Theta}=\{\tau(\mu)|\ \mu\in Q(\pi)$ and 
$\Theta(\mu)=\mu\}$.  The involution $\theta$ of $\gg$ lifts to 
  a
  ${\bf C}$ algebra isomorphism $\tilde\theta$ of $U$ ([L3, Theorem 
  7.1]). 
  Set $$B_i=y_it_i+\tilde\theta(y_i)t_i$$ for all $\alpha_i\in 
  \pi\setminus \pi_{\Theta}$.   The standard quantum analog of 
  $U(\gg^{\theta})$ inside of $U_q(\gg)$ is the  algebra 
  $B_{\theta}$ generated by ${\cal M},T_{\Theta}$, and $B_i$ for 
  $\alpha_i\in \pi\setminus\pi_{\Theta}$.

  For most irreducible symmetric pairs, any quantum analog of   
  $U(\gg^{\theta})$ inside of $U_q(\gg)$ is isomorphic to 
  $B_{\theta}$ via a Hopf algebra automorphism of $U_q(\gg)$.   
  However, under certain circumstances, $U(\gg^{\theta})$ admits a 
  one parameter family of analogs.  
     Define 
two  subsets ${\cal S}$ and ${\cal D}$ of $\pi^*$ by 
$${\cal S}=\{\alpha_i\in \pi^*| \ \Theta(\alpha_i)=-\alpha_i,{\rm \ 
and \ if\ }\Theta(\alpha_j)=-\alpha_j {\rm \ then \ }
{{(\alpha_i,\alpha_j)}\over{(\alpha_j,\alpha_j)}}\in {\bf Z}\}$$ and
 $${\cal D}=\{\alpha_i\in \pi^*|\ i\neq {\rm p}(i),{\rm \ and\ }(\alpha_i,\Theta(\alpha_i))\neq 0\}.$$
Since we are assuming that $\gg,\gg^{\theta}$ is irreducible, it 
follows that ${\cal S}\cup{\cal D}$ has at most one element ([L4, 
Section 2 and Section 7]).
Suppose first that  ${\cal S}$ is not empty.  Given $s\in {\cal C}$,
  let $B_{\theta,s,1}$ 
be the subalgebra of $U_q(\gg)$ generated by ${\cal 
M}$, $T_{\Theta}$,  and $B_{i}$ for $\alpha_i\in 
\pi\setminus (\pi_{\Theta}\cup {\cal S})$  and $B_{j,s,1}$ for 
$\alpha_j\in {\cal S}$  where  
$$B_{j,s,1 }= 
y_jt_j+\tilde\theta(y_j)t_j +st_j.$$   Note that when $s=0$, 
  the 
algebra $B_{\theta,0,1}$ is just the standard quantum analog 
$B_{\theta}$.  Now consider the 
case when ${\cal D}$ is nonempty.   Then given $d\in {\cal C}$ with 
$d\neq 0$, let $B_{\theta,0,d}$ 
be the subalgebra of $U_q(\gg)$ generated by ${\cal 
M}$, $T_{\Theta}$,   $B_{i}$ for $\alpha_i\in 
\pi\setminus (\pi_{\Theta}\cup {\cal D})$  and $B_{j,0,d}$ for 
$\alpha_j\in {\cal D}$  where  
$$B_{j,0,d }= 
y_jt_j+d\tilde\theta(y_j)t_j .$$  In the special case when $d=1$, we 
have that $B_{\theta,0,1}$ is just the standard analog $B_{\theta}$.
Given $\alpha_i\in {\cal S}$,  $B_{i,s,1}$ 
is shortened 
to $B_i$ when $ { s}$  can be read from context.  A similar 
convention is applied to $B_{i,0,d}$ for $\alpha_i\in {\cal D}$.  Set 
$s_i=s$ for $\alpha_i\in {\cal S}$ and $s_i=0$ for all 
$\alpha_i\notin{\cal S}$. Similarly, set $d_i=d$ for $\alpha_i\in 
{\cal D}$ and  set $d_i=1$ for all $\alpha_i\notin{\cal 
D}$.    In order to make arguments in the general case, it will often be 
convenient to write each $B_i$ as 
$$B_i=y_it_i+d_i\tilde\theta(y_i)+s_it_i.$$

Let $U_{\cal R}$ be the ${\cal R}$ subalgebra of $U$ generated by 
$x_i,y_i,t_i^{\pm 1}$ for $1\leq i\leq n$.    A subalgebra $B$ of 
  $U_q(\gg)$ is called real if $B=(B\cap 
U_{\cal R})\oplus i(B\cap 
U_{\cal R})$.  Note that the standard analog $B_{\theta}$ is real 
since it is generated by elements of $U_{\cal R}$.  
More generally, when ${\cal S}$ is nonempty, $B_{\theta,s,1}$ is real 
provided $s\in {\cal R}$.   Similarly, when ${\cal D}$ is nonempty, 
$B_{\theta,0,d}$ is real provided $d\in {\cal R}$.

We define a family ${\cal B}$ of coideal subalgebras of $U$ associated to 
$\gg,\gg^{\theta}$ as follows.   Let ${\bf H}$ denote the set of Hopf algebra automorphisms of $U$ which fix 
elements of $T$.  Note that ${\bf H}$ acts on the set of coideal 
subalgebras of $U$. If ${\cal S}\cup {\cal D}$ is empty, set ${\cal B}$ 
equal to the orbit of $B_{\theta}$ under ${\bf H}$. If ${\cal S}$ is 
nonempty, set ${\cal B}$ equal to the union of the orbits of the 
set $\{B_{\theta,s,1}|s\in {\cal R}\}$ under ${\bf H}$.  If ${\cal D}$ is 
nonempty, set ${\cal B}$ equal to the union of the orbits of the 
set $\{B_{\theta,0,d}|d\in {\cal R}$ and $d\neq 0\}$ under ${\bf H}$.     
By [L3, Theorem 7.5],    any real quantum 
analog of $U(\gg^{\theta})$ inside $U$ is isomorphic via a Hopf 
algebra automorphism in ${\bf H}$ to some element 
of ${\cal B}$.   It should be noted, however, that  not every element in 
${\cal B}$ is a quantum analog of $U(\gg^{\theta})$ inside $U$ since 
not all the subalgebras in ${\cal B}$ specialize to $U(\gg^{\theta})$ as $q$ goes 
to $1$.  For example, if ${\cal S}$ is nonempty, then 
$B_{\theta,s,1}$ specializes to $U(\gg^{\theta})$ if and only if $s$ 
specializes to $0$.   Despite this fact, we will refer to the set ${\cal B}$ as the 
orbit under ${\bf H}$ of  the set of real 
quantum analogs of $U(\gg^{\theta})$ inside of $U(\gg)$.

 Let $F_r(\check U)$ denote the locally finite part of $\check U$ 
 with respect to the right adjoint action ${\rm ad}_r$.    It follows from 
 [L4, 
 Theorem 3.1], that for each $B\in {\cal B}$ there exists a   conjugate linear  algebra 
 antiautomorphism of $U$ which restricts to an antiautomorpism of 
 $B$.  Thus we have the following version of [L3, Theorem 7.6] and 
 [L2, Theorem 3.5].  

 \begin{theorem}   Given $B\in {\cal B}$, the action of $B$ on   $F_r(\check U)$
    is semisimple.   Moreover, if $V$ is a finite dimensional $(\adr 
     B)$ submodule of $\check U$, then $V$ is a subset of $F_r(\check 
     U)$.  
 \end{theorem}
 
For each $B\in {\cal B}$, let $\check U^B$ denote the set of $B$ 
invariants of $\check U$ with respect to the right adjoint action.  
By [L5, Lemma 3.5], $\check U^B$ is just the centralizer in $\check U$ of $B$.   
It follows that the center  $Z(\check U)$ of $\check U$ is contained 
in $\check U^B$. Since $\check U^B$ is just a direct sum of one-dimensional simple $B$ 
modules, we see  from the above theorem that $\check U^B$ is a submodule of $F_r(\check U)$.

  For each $\lambda\in P(\pi)\cup P(\Sigma)$, let $z^{\lambda}$ be the algebra 
  homomorphism from  
  $\check U^0$ to ${\cal C}$  defined by $z^{ \lambda}(\tau(\gamma))=
  q^{(\lambda,\gamma)}$ for all $\gamma\in P(\pi)$. 
  Given 
  $\Lambda\in {\rm Hom}(\check U^0,{\cal C})$, let $L(\Lambda)$ denote the 
  simple $U$ module with highest weight $\Lambda$.  If 
  $\Lambda=z^{\lambda}$ for some $\lambda\in P(\pi)$, then we write 
  $L(\Lambda)$ as $L(\lambda)$.
  
  For any multiplicative group $G$, we write ${\cal C}[G]$ for the 
  group algebra  generated by $G$ over ${\cal C}$.   This notation 
  will be applied to groups as well as multiplicative monoids related 
  to
  $\check T$.
  Now suppose that $H$ is an additive group such as $P(\pi)$ or 
  $Q(\Sigma)$.  Let
  $<z^{\lambda}|\ \lambda\in H>$ denote the multiplicative group isomorphic to 
  $H$ via the map which sends  $z^{\lambda}$ to $\lambda$.   
We write ${\cal C}[H]$ for   the group algebra of the multiplicative 
  group $<z^{\lambda}|\ \lambda\in H>$.
  
  An appendix with the definition of symbols introduced after this 
  section can be found at the end of the paper.   For undefined 
  notions or more information about notation, the reader is referred 
  to [Jo] or [L5].

  \section{ The Harish-Chandra Map}

In this section, we construct a quantum 
Harish-Chandra map   associated  to the symmetric pair $\gg,\gg^{\theta}$.
The procedure involves  four  tensor product decompositions 
of the quantized enveloping algebra with respect to a subalgebra 
$B\in {\cal B}$.  The results establishing these decompositions 
parallels the material presented in [L5, Section 2].  However, the 
theorems in [L5, Section 2] are only proved for  
  the standard analogs of 
 $U(\gg^{\theta})$ which have a reduced  restricted root system 
$\Sigma$. In particular, we 
first generalize [L5, Lemma 2.1] and [L5, Theorem 2.2] so that it 
applies to all coideal subalgebras in ${\cal B}$ associated to any 
irreducible symmetric pair $\gg,\gg^{\theta}$.  The quantum 
Harish-Chandra map  associated to an algebra $B\in {\cal B}$  is then defined using these decompositions.  The 
end of the section focuses on the restriction of this map to the 
$(\adr B)$ invariants of $\check U$.
The tensor product  decompositions of this section   
 are also   used in Section 4 to 
generalize  the  notion of  radial components  defined in [L5, 
Section 3].

We introduce some notation from    [L5, Section 
2]. Set ${\cal M}^+={\cal M}\cap U^+$ and ${\cal M}^-={\cal M}\cap 
G^-$.  Write ad for the left adjoint action of $U$ on itself. 
Let $N^+$ be the subalgebra of $U^+$
generated by the $(\ad {\cal M}^+)$ module $(\ad {\cal M}^+){\cal
C}[x_i |\alpha_i\notin \pi_{\Theta}].$   Similarly, let $N^-$ be the subalgebra
  of $G^-$ generated by the $(\ad {\cal
M}^-)$ module $(\ad {\cal M}^-){\cal C}[y_it_i |\alpha_i\notin
\pi_{\Theta}]$. 
Given $\beta\in Q(\pi)$ and a subspace $S$ of $U$, we write $S_{\beta, 
r}$ for the sum of the  weight subspaces $S_{\beta'}$ of $S$ with 
$\tilde\beta'=\tilde\beta$.   
  Set $T'_{\geq}$ equal to the multiplicative monoid generated by the 
$t_i^2$ for $\alpha_i\in \pi^*\setminus {\cal S}$ and $t_i$ for 
$\alpha_i\in {\cal S}$.  Note that [L4, Lemma 3.5] implies that $T'_{\geq}$ is just equal to the 
set $\{\tau(\gamma)|\ \gamma\in Q^+(\pi) $ and $\tilde\gamma\in   P(2\Sigma)\}$.
 We have the following generalizations of [L5, Lemma 2.1] and [L5, Theorem 2.2].

\begin{lemma}
For each $B\in {\cal B}$, all $\beta, \gamma\in Q^+(\pi)$,
and $Y\in U^+_{\gamma}G^-_{-\beta}$, we have
\begin{enumerate}
    \item[(i)]  $Y\in
N^+_{\beta+\gamma,r}B+ \sum_{\tilde\beta'<\tilde\beta+\tilde\gamma}
N^+_{\beta',r} T'_{\geq}B  $
\item[(ii)]
 $Y\in
BN^-_{-\beta-\gamma,r}+\sum_{\tilde\beta'<\tilde\beta+\tilde\gamma} B T'_{\geq}N^-_{-\beta',r}.$
\end{enumerate}
\end{lemma}

\noindent
{\bf Proof:}  Fix $B\in {\cal B}$.
We consider here the proof of the first assertion; the second 
assertion follows with a similar argument.   As in [L5], we may 
reduce to the case where $Y$ is a monomial of the form  
$y_{i_1}t_{i_1}\cdots y_{i_m}t_{i_m}$  such 
that   $y_{i_m}t_{i_m}$ is not   in 
$B$. Set $\beta$ equal to the weight of $Y$.   The proof in [L5] proceeds by 
induction on $m$.       In 
particular, we assume that all monomials in the $y_it_i$ of length 
less than or equal to $m-1$ satisfy (i).  
 Note that 
$$\eqalign{Y=
&y_{i_1}t_{i_1}\cdots y_{i_{m-1}}t_{i_{m-1}}
B_{i_m} 
 -d_{i_m}y_{i_1}t_{i_1}\cdots
y_{i_{m-1}}t_{i_{m-1}}\tilde\theta(y_{i_m})t_{i_m}\cr
&-s_{i_m}y_{i_1}t_{i_1}\cdots y_{i_{m-1}}t_{i_{m-1}}t_{i_m}.\cr}$$
By the proof of [L5, Lemma 2.1], 
both $$y_{i_1}t_{i_1}\cdots y_{i_{m-1}}t_{i_{m-1}}
B_{i_m}{\rm\ and \ }y_{i_1}t_{i_1}\cdots
y_{i_{m-1}}t_{i_{m-1}}\tilde\theta(y_{i_m})t_{i_m}$$ are contained in 
the right hand side of (i). 

Now suppose that $s_{i_m}\neq 0$.  It follows  that $\alpha_{i_m}\in {\cal S}$ and 
thus $t_{i_m}\in T'_{\geq}$.
Hence $$s_{i_m}y_{i_1}t_{i_1}\cdots y_{i_{m-1}}t_{i_{m-1}}t_{i_m}\in 
{\cal C}[T'_{\geq}]y_{i_1}t_{i_1}\cdots y_{i_{m-1}}t_{i_{m-1}}.$$  
Note that $N_{\beta',r}^{+}T'_{\geq}=T'_{\geq}N^{+}_{\beta',r}$
and $N_{\beta',r}^{-}t=tN^{-}_{\beta',r}$ for each 
$\beta'$ and for each $t\in T'_{\geq}$. The lemma now follows by applying the inductive hypothesis 
to $y_{i_1}t_{i_1}\cdots y_{i_{m-1}}t_{i_{m-1}}$. $\Box$ 

\medskip
Let $T'$ be the subgroup of $T$ generated by $t_i$ for $\alpha_i\in 
\pi^*$.  
The next result provides four tensor product decompositions of $U$.
This generalizes [L5, Theorem 2.2] which is proved for the subset 
of standard analogs in ${\cal B}$.  Though the proof here is quite 
similar to the argument in [L5], one of the antiautomorphisms used in 
[L5]  does not work in 
this more general setting.

 \begin{theorem} For all $B\in {\cal B}$, there
are vector space isomorphisms via the multiplication map
\begin{enumerate}
\item [(i)]$   N^+ \otimes {\cal C}[T']\otimes B \cong U $ 
\item [(ii)] $B\otimes {\cal C}[T']\otimes N^+\cong U$
\item [(iii)]$ N^- \otimes {\cal C}[T']\otimes B \cong U $
\item [(iv)] $ B\otimes {\cal C}[T']\otimes N^-\cong U$
\end{enumerate}
\end{theorem}

\noindent
{\bf Proof:}   Fix $B\in {\cal B}$. 
By [L3, Section 6],  the multiplication map induces an isomorphism 
$$U^+\cong {\cal M}^+\otimes N^+.\eqno{(2.1)}$$  
 Switching the order of ${\cal 
M}^+$ and $N^+$, the same argument  shows that 
multiplication  induces an isomorphism 
$$U^+\cong  N^+\otimes {\cal M}^+.\eqno{(2.2)}$$ 
The triangular decomposition of $U$ ([L5, (2.6)]) and (2.1)
imply that multiplication induces an isomorphism (see [L5, (2.8)]) $$U\cong  G^-\otimes   {\cal
M}^+\otimes {\cal C}[T_{\Theta}]\otimes{\cal C}[T'] \otimes
N^+.$$  A similar argument   using (2.2)  yields that 
multiplication also induces an isomorphism 
$$U\cong  N^+\otimes  {\cal C}[T']\otimes {\cal C}[T_{\Theta}] 
\otimes{\cal M}^+ \otimes
G^-.$$ 
As explained in [L5] (see [L5,(2.9)] and subsequent discussion), the lowest weight 
term of  an element $b\in B$ written as a 
sum of weight vectors  is in $G^-{\cal M}^+T_{\Theta}$.  By [L3, 
Section 6 and (6.2)] and [Ke], we have $G^-{\cal 
M}^+T_{\Theta}=N^-{\cal M}T_{\Theta}={\cal M}T_{\Theta}N^-={\cal M}^+T_{\Theta}G^-$. Hence 
$Bv\cap Bv'=0$ and $vB\cap v'B=0$ for any two linearly independent 
elements $v$ and $v'$  of $T'N^+$.  Therefore the multiplicaton map 
induces injections from $B\otimes {\cal C}[T']\otimes N^+$ and 
$N^+\otimes {\cal C}[T']\otimes B$ into $U$.
  
Recall that there is a antiautomorphism $\kappa$ of $U$ such 
that $\kappa(N^-)=N^+$, $\kappa(B)=B$, and $\kappa(U^+)=G^-$
(see [L4, Theorem 3.1] and the proof of [L5, Theorem 2.2]).  
Applying $\kappa$ to the above yields injections from 
$B\otimes {\cal C}[T']\otimes N^-$ and 
$N^-\otimes {\cal C}[T']\otimes B$ into $U$.

By Lemma 2.1(i), $U^+G^-\subseteq N^+T'B$.   Since $T=T'\times 
T_{\Theta}$, it follows that $U=U^+G^-T\subseteq N^+T'B$.  Hence
$U=N^+T'B$  and isomorphism (i) now follows.
Assertion  (iii) follows similarly from Lemma 2.1(ii).
Isomorphism (ii) now follows from (iii) and (iv) from (i) using the
antiautomorphism $\kappa$ as in the proof of 
[L5,
Theorem 2.2].
$\Box$

\medskip

We introduce two important groups related to $\check T$ that will be 
used throughout the paper.
Note that  the group $\check T$ 
can be extended to a larger group isomorphic to the additive group 
$\{\mu/2|\ \mu\in 
P(\pi)\}$ via the obvious extension of $\tau$.  Let $\check {\cal A}$ 
and $\check T_{\Theta}$ 
be the subgroups of this extension defined by 
$\check {\cal A}=\{\tau(\tilde\mu)|\mu\in
P(\pi)\} $ 
and    $\check T_{\Theta}=\{\tau(1/2(\mu+\Theta(\mu))|\ \mu\in 
P(\pi)\}.$    Consider for the moment
$a\in \check U^B$.  It follow that $a$ is a sum of weight vectors 
each of whose weight $\beta$ satisfies $\Theta(\beta)=-\beta$.  Thus the 
invariant ring $\check U^B$ is equal to $\check U^{B\check T_{\Theta}}$.

We are now ready to define the quantum Harish-Chandra projection map 
with respect to a subalgebra $B$ in ${\cal B}$.
Note that $\tau(\mu)=\tau(\tilde\mu)\tau(1/2(\mu+\Theta(\mu))\in 
\check {\cal A}\check T_{\Theta}$ for all $\mu\in P(\pi)$. 
Hence we have  the inclusion:
$$\check U^0\subset {\cal C}[\check{\cal A}]\oplus
\check U^0{\cal C}[\check T_{\Theta}]_+.\eqno{(2.3)}$$
Theorem 2.2 and (2.3)   imply the following inclusion for each 
$B\in {\cal B}$ 
$$\check U\subseteq ((B\check T_{\Theta})_+\check U+N^+_+\check {\cal A})\oplus 
{\cal C}[\check A].\eqno{(2.4)}$$

\begin{definition}The (quantum) Harish-Chandra map with respect to 
the symmetric pair $\gg,\gg^{\theta}$ and subalgebra $B$ in ${\cal B}$ is 
the projection   ${\cal P}_B$  of $\check U$ onto ${\cal C}[\check{\cal 
A}]$ using the direct sum  
decomposition (2.4).
\end{definition}

 It should be noted that the map ${\cal P}_B$ is the same as the map 
 ${\cal P}_{\cal A}$ defined for    standard analogs of $U(\gg^{\theta})$ where $\Sigma$ 
 is a reduced root system in [L5, immediately following (3.11)].   
Note further that   ${\cal P}_B$ is a linear map for each $B\in {\cal B}$.
However, ${\cal P}_B$ is not an algebra homomorphism on $\check U$.   
The next result shows that ${\cal P}_B$ is an algebra homomorphism upon
restriction to $\check U^B$.

\begin{theorem}  For all $B\in {\cal B}$ the restriction of ${\cal P}_{B}$ to 
$\check U^B$ is an algebra homomorphism.

\end{theorem}   

\noindent{Proof:}  Suppose that  $a$ and $a'$ in $\check U^B$.    
We have $$aa'\in a((B\check T_{\Theta})_+\check U+N_+^+\check {\cal A})+a{\cal P}_B(a').$$  Since $a\in \check U^B$, it follows 
that $a(B\check T_{\Theta})_+\check U\subseteq 
(B\check T_{\Theta})_+\check U$. Using  (2.4) we have $$aN^+_+\check {\cal A}\subset
\check UN^+_+\check {\cal A}=(B\check 
T_{\Theta})_+\check U+N^+\check{\cal A})N^+_+\check {\cal A}\subseteq ((B\check 
T_{\Theta})_+\check U+N^+_+\check {\cal A}.$$ Hence 
$$\eqalign{aa'&\in (B\check T_{\Theta})_+\check U+(aN_+^+\check{\cal A})+a{\cal 
P}_B(a')\cr
&\subseteq (B\check T_{\Theta})_+\check U+ N_+^+\check{\cal A})+{\cal P}_B(a){\cal 
P}_B(a').\cr}$$  It follows that ${\cal P}_B(aa')={\cal P}_B(a){\cal 
P}_B(a')$.   The theorem now follows from the fact that 
  ${\cal P}_B$ is a linear map. $\Box$
  
\medskip Let ${\cal A}$ denote the subgroup of  $\check 
{\cal A}$  consisting of those elements 
$\tau(2\mu)$ with $\mu\in P(\Sigma)$.   (Note that this definition 
of ${\cal A}$ differs from the definition of ${\cal A}$ in [L5].)  
Let $G^+$ denote the subalgebra of $U$ generated by $x_it_i^{-1}$ for 
$1\leq i\leq n$ and 
let $U^-$ denote the subalgebra of $U$ generated by $y_i$  for $1\leq 
i\leq n$.  

\begin{lemma}  Suppose that $u\in \check U^B$.   Then $$u\in 
 U^+G^-{\cal A}+U^+G^-{\cal C}[\check T_{\Theta}]_+.$$
 \end{lemma}
 
 \noindent
 {\bf Proof:}   Recall (Theorem 1.1) that $\check U^B$ is a subset of $F_r(\check 
 U)$.   Moreover,   $F_r(\check U)$ is a direct sum of $\adr U$ 
 submodules of the form $(\adr U)\tau(2\gamma)$ with $\gamma\in 
 P^+(2\pi)$ (see [Jo, Section 7] and [L5, discussion preceding Lemma 7.2]).     Thus it is sufficient to prove the lemma 
 for $u\in \check U^B\cap (\adr U)\tau(2\gamma)$ where 
 $\gamma\in P^+(\pi)$.  As explained in [L5,  discussion 
 preceding Lemma 7.2], we have that 
 $$(\adr U)\tau(2\gamma)\subset \sum_{\mu\in 
 Q^+(\pi)}G^+U^-\tau(2\gamma-2\mu).\eqno{(2.5)}$$
 
 Now $(\adr U)\tau(2\gamma)$ can be written as a direct sum of 
 simple $(\adr U)$ modules.   In particular, $u$ is a sum of $(\adr 
 B)$ invariant vectors contained in  
  simple $(\adr U)$ submodules of $(\adr 
 U)\tau(2\gamma)$. It follows from [L4, Theorem 3.6] that   $u=\sum_{\tilde\lambda\in 
 P(\Sigma)}w_{2\tilde\lambda}$ for some weight vectors  $w_{2\tilde\lambda}$  in 
 $(\adr U)\tau(2\gamma)$ of weight $2\tilde\lambda$. Using (2.5), we can write
 $$w_{2\tilde\lambda}\in \sum_{\mu\in Q^+(\pi)}\sum_{\alpha- 
 \beta=2\tilde\lambda}G^+_{\alpha}U^-_{-\beta}\tau(2\gamma-2\mu)$$ for each 
  $\tilde\lambda\in P(\Sigma)$.  Now consider $\alpha$ and $\beta$ in 
 $Q^+(\pi)$ such that $\alpha-\beta=2\tilde\lambda$.  Note that
 $G^+_{\alpha}=U^+_{\alpha}\tau(- \alpha)$ and 
 $U^-_{\beta}=G^-\tau(-\beta)$.  Hence 
 $G^+_{\alpha}U^-_{\beta}=U^+G^-\tau(-2\alpha)\tau(\alpha-\beta).$
 Since $\alpha-\beta=2\tilde\lambda$, it follows that 
 $\tau(\alpha-\beta)\in {\cal A}$.   Therefore, 
 $$u\in \sum_{\eta\in Q(\pi)}U^+G^-\tau(2\eta){\cal A}.$$
 Now $\tilde\eta\in Q(\Sigma)$ and hence $\tau(2\tilde\eta)\in {\cal A}$ for all $\eta\in Q(\pi)$.  
     The lemma now follows from the decomposition  in (2.3). 
 $\Box$
 \medskip

 Recall the definition of $T'_{\geq}$ given before Lemma 2.1.
 By [L4, Lemma 3.5], $\tilde\alpha_i/2\in P(\Sigma)$ whenever 
 $\alpha_i\in {\cal S}$. It follows that $\tau(\tilde\mu)\in {\cal 
 A}$ whenever $\tau(\mu)\in T'_{\geq}$.   Hence by (2.3) and the 
 discussion preceding it, we have that ${\cal P}_B(T'_{\geq})\subseteq 
 {\cal A}$.  This fact combined with the previous lemma allows us 
 to refine   the codomain of ${\cal P}_B$.
 
 \begin{theorem} For all $B\in {\cal B}$, the image of $\check U^B$ 
 under ${\cal P}_B$ is contained in ${\cal C}[{\cal A}]$.
 \end{theorem}
 
 \noindent
 {\bf Proof:} Let $u\in \check U^B$.   By Lemma 2.5, 
 there exists $u'\in  U^+G^-{\cal A}$ such that ${\cal P}_B(u)={\cal 
 P}_B(u')$. Note that $U^+G^-{\cal A}={\cal A}U^+G^-$. It follows from 
 Lemma 2.1(i) and the above discussion that ${\cal P}_B(u')\in {\cal C}[{\cal A}]$. $\Box$

\section{A comparison of two root systems}

Let $B\in {\cal B}$ and consider for the moment the restriction of ${\cal P}_B$ to $\check 
U^0$.   The image of an element $\tau(\mu)$   under 
${\cal P}_B$
  is just $\tau(\tilde\mu)$. Thus there is a    close 
connection between the quantum Harish-Chandra map ${\cal P}_B$ and the 
map $\ \tilde{\ }\ $ on $\gh^*$.    This section is an analysis 
of the latter map.
The image of the 
fundamental weights associated to $\pi$ under the restriction map 
$\ \tilde{\ }\ $ are completely determined here.  As a consequence, we 
   establish necessary and 
sufficient conditions  for $\widetilde{P^+(\pi)}$ to equal 
$P^+(\Sigma)$.  This result is  used  to understand the image of  
 the center $Z(\check U)$ under  ${\cal P}_{B}$ in 
Sections 5 and 6.  

For each $1\leq i\leq n$, set $\omega_i$ equal to the fundamental weight corresponding to the 
simple root $\alpha_i$. Recall that the restricted root system is 
either an ordinary reduced root system as classified in [H, Chapter 
III] or is 
nonreduced of type 
BC$_r$ for some integer $r\geq 1$ ([Kn, Chapter II, Section 5]).   In the latter case,   there is a 
unique simple root $\alpha_j\in \pi^*$ such that 
both $\tilde\alpha_j$ and 
$2\tilde\alpha_j$ are restricted roots.         Let
$\omega'_i$ denote the fundamental weight corresponding to the simple
root $\tilde\alpha_i$  with respect to the restricted root system.  
In particular, if $2\tilde\alpha_i$ is not a restricted root then 
$(\omega_i',\tilde\alpha_i)=(\tilde\alpha_i,\tilde\alpha_i)/2$ and 
if $2\tilde\alpha_i$ is a restricted root then 
$(\omega_i',2\tilde\alpha_i)=(2\tilde\alpha_i,2\tilde\alpha_i)/2$.

   We can break up the set $\{\tilde\alpha_i|\ \alpha_i\in \pi^*\}$ 
   into three cases.  In each case, we give the value of  
   $(\tilde\alpha_i,\tilde\alpha_i)$.
    
   \medskip
   \noindent
    {\bf Case 1:} $\Theta(\alpha_i)=-\alpha_i$.   Then 
    $(\tilde\alpha_i,\tilde\alpha_i)=(\alpha_i,\alpha_i)$.
    
    \medskip
    \noindent
    {\bf Case 2:} $\Theta(\alpha_i)\neq -\alpha_i$, 
    and $(\alpha_i,\Theta(\alpha_i))=0$.  Then 
    $\tilde\alpha_i=(\alpha_i-\Theta(\alpha_i))/2$ and 
    $(\tilde\alpha_i,\tilde\alpha_i)=(\alpha_i,\alpha_i)/2$.
    
    \medskip
    \noindent
    {\bf Case 3}: $\Theta(\alpha_i)\neq -\alpha_i$, 
    and $(\alpha_i,\Theta(\alpha_i))\neq 0$.  In this case,
    $\alpha_i+\Theta(-\alpha_i)$ is a 
    root. Note that $\tilde\alpha_i=\widetilde{\Theta(-\alpha_i)}$. 
    Hence, both 
    $\tilde\alpha_i$ and $2\tilde\alpha_i$ are roots in $\Sigma$.   
    In particular, as an abstract root system, $\Sigma$ must 
    be of type BC$_r$ for some $r$.   The only way this can happen is 
    for 
    $(\alpha_i,\Theta(\alpha_i))=-(\alpha_i,\alpha_i)/2$.  It follows 
    that $(\tilde\alpha_i,\tilde\alpha_i)=(\alpha_i,\alpha_i)/4$.

\begin{lemma}  Supppose $\alpha_i\in \pi^*$.  If $i= {\rm p}(i)$ and 
$\alpha_i\neq \Theta(-\alpha_i)$ then $\tilde\omega_i=2\omega_i'$.   
If $i\neq {\rm p}(i)$ or if $\alpha_i=\Theta(-\alpha_i)$ then $\tilde\omega_i=\omega_i'$.  
\end{lemma}

\noindent
{\bf Proof:}
Fix $\alpha_i\in \pi^*$.   Note that 
$(\tilde\omega_i,\tilde\alpha_j)=(\omega_i,\tilde\alpha_j)=0$ 
whenever $\tilde\alpha_j\neq \tilde\alpha_i$.   It follows that 
$\tilde\omega_i$ is a scalar multiple of $\omega_i'$.  Now   
$$(\tilde\omega_i,\tilde\alpha_i)=(\omega_i,\tilde\alpha_i)=(1+\delta_{i{\rm p}(i)})(\alpha_i,\alpha_i)/4.$$
Note that we can break up both Case 2 and Case 3 into two subcases 
depending on whether $i={\rm p}(i)$ or $i\neq {\rm p}(i)$. The lemma follows from a straightforward computation using the three 
cases above and the corresponding subcases.
$\Box$

\medskip

For each $i$ such that $\alpha_i\in \pi^*$, set $\pi_i=\{\alpha_j|\ 
(\omega_j,\tilde\alpha_i)\neq 0\}$.  Let $\gg_i $ be the semisimple 
Lie subalgebra of $\gg$ generated by the positive and negative root 
vectors 
corresponding to roots in $\pi_i$.   Note that $\theta$ restricts to 
an involution of $\gg_i$ and that $\gg_i,\gg_i^{\theta}$ is   a rank 
one symmetric pair.  In particular, 
the restricted root system $\Sigma_i$ associated to 
$\gg_i,\gg_i^{\theta}$ has rank one.  Moreover 
$\Sigma_i=\{\pm\tilde\alpha_i\}$ if $2\tilde\alpha_i$ is not a 
restricted root and 
$\Sigma_i=\{\pm\tilde\alpha_i,\pm2\tilde\alpha_i\}$  
otherwise.
Let $\beta_i$ denote the longest positive root in $\Sigma_i$.  Thus, if $\Sigma_i$ is reduced then $\beta_i=\tilde\alpha_i$ 
and if $\Sigma_i$ is not reduced then $\beta_i=2\tilde\alpha_i$.

\begin{lemma} 
    For each $\alpha_i\in \pi^*$, there  exists 
    $\lambda_i\in P(\pi)$ such that $(\lambda_i,\beta_j)=(\beta_j,\beta_j)/2$.
    \end{lemma}

\noindent 
{\bf Proof:}  
The proof of this lemma is a straightforward computation. The list 
below of possibilities for the rank one symmetric pair 
$\gg_i,\gg_i^{\theta}$  can be deduced from 
  Araki's classification ([A], see also [L5, Section 4]) of irreducible symmetric 
  pairs. For cases (3.3) through (3.7), we relabel the set $\pi_i$ as 
  $\{\alpha_{i1},\dots,\alpha_{ir}\}$ using the second subscript to 
  indicate the  numbering of the simple roots for a particular 
  root system as given in [H].  Given $\alpha_{ij}\in \pi_i$, we 
  write $\omega_{ij}$ for the fundamental weight in $P(\pi)$ 
  corresponding to the simple root $\alpha_{ij}$ considered as an 
  element in the larger set $\pi$.
  For each $i$ with $\tilde\alpha_i\in 
\pi_{\Theta}$, we
provide at least one choice (which might not be the only choice) of  
$\lambda_i$ in $P(\pi)$ satisfying the 
desired
condition.  

\begin{tabbing}
	\medskip
\noindent
\=$\gg_i$ is of \=type $A_1$ with $\pi_i=\{\alpha_i\}$,
$\Theta(\alpha_i)=-\alpha_i$, and $\lambda_i=\omega_i$\qquad\qquad\quad 
\=(3.1)\\

\> $\gg_i$ is of \>type $A_1\times A_1$ with 
$\pi_i=\{\alpha_i,\alpha_{{\rm p}(i)}\}$,
$\Theta(\alpha_i)=-\alpha_{{\rm p}(i)}$, and\>(3.2)\\
\>\>$\lambda_i   =\omega_i$.\\ 

\>$\gg_i$ is of \>type $A_3$ with
$\pi_i=\{\alpha_{i1},\alpha_{i2},\alpha_{i3}\}$,
$\alpha_i=\alpha_{i2}$,\>(3.3)\\
\>\>$\Theta(\alpha_{i2})=-\alpha_{i1}-\alpha_{i3}-\alpha_{i2}$,
and $\lambda_i= \omega_{i1}$ or $\lambda_i=\omega_{i3}$. \>\\ 

\>$\gg_i$ is of \>type $A_r$ with $r\geq 2$, $\pi_i=\{\alpha_{i1},\dots, 
\alpha_{ir}\}$,
$\alpha_i=\alpha_{i1}$, \>(3.4)\\
\>\> $\Theta(\alpha_{ij})=-\alpha_{i,j'}-\alpha_{i2}-\dots-\alpha_{i,r-1}$ 
where $\{j,j'\}= \{1,r\}$,\\
\>\>and $\lambda_i= \omega_{i1}$.\>\\

\>$\gg_i$ is of \>type $B_r$ with
$\pi_i=\{\alpha_{i1},\dots,\alpha_{ir}\}$, $\alpha_i=\alpha_{i1}$, 
\>(3.5)\\
\>\>$\Theta(\alpha_{i1})=- \alpha_{i1}-2\alpha_{i2}-\cdots 
-2\alpha_{ir}$, and $\lambda_i= \omega_{ir}$.\>\\

\>$\gg_i$ is of \>type $D_r$ with
$\pi_i=\{\alpha_{i1},\dots,\alpha_{ir}\}$, $\alpha_i=\alpha_{i1}$,\>(3.6)\\
\>\>$\Theta(\alpha_{i1})=-\alpha_{i1}-2\alpha_{i2}-\cdots
-2\alpha_{i,r-2}-\alpha_{i,r-1}-\alpha_{i,r}$,\\ \>\>and 
 $\lambda_i= \omega_{ir}$. \>\\

\>$\gg_i$ is of \>type $C_r$ with $\pi_i=\{\alpha_{i1},\dots, 
\alpha_{ir}\}$,
$\alpha_i=\alpha_{i2}$,  \>(3.7)\\
\>\>$\Theta(\alpha_i)=-\alpha_{i1}- \alpha_{i2}-2\alpha_{i3}\dots 
-2\alpha_{i,r-1}-\alpha_{ir}$, and $\lambda_i= \omega_{ir}$.\>\\

\>$\gg_i$ is of \>type $F4$ with $\gg=\gg_i$, $\alpha_i=\alpha_4$,\>(3.8)\\ 
\>\>$\Theta(\alpha_4)=-\alpha_4-3\alpha_3-2\alpha_2-\alpha_1$ 
 and $\lambda_i= \omega_4$.$\Box$\>\\
\end{tabbing}

 Let $W$ denote the Weyl group associated to the root system
of $\gg$  and let $W_{\Theta}$ denote the Weyl group associated to 
the restricted root system $\Sigma$.  It is well known that 
the restricted Weyl group $W_{\Theta}$ is a homomorphic image of the subgroup of W which leaves
$Q(\Sigma)$ invariant.    
 (See   [He, p. 791].)  
In particular,
given $\tilde w$ in $W_{\Theta}$, there exists $ w $ in $W$ such that  
$\tilde w$ and $w$ act the same on $Q(\Sigma)$.

In the next two lemmas, we examine the relationship between the 
integral weights associated to $\pi$ and those associated to 
$\Sigma$.

\begin{lemma}  The set $\widetilde {P(\pi)}$   is a 
subset of $P(\Sigma)$.
    \end{lemma}
    
    \noindent{Proof:}   We need to show that each $\beta\in 
    P(\pi)$ satisfies $(\tilde\beta,\tilde\alpha)\in {\bf 
    Z}(\tilde\alpha,\tilde\alpha)/2$ for all $\tilde\alpha\in \Sigma$.  Now suppose $\tilde\alpha\in 
    \Sigma$.   We can find $\tilde w\in W_{\Theta}$ such that 
    $\tilde w\tilde\alpha=\tilde\alpha_i$ or $\tilde 
    w\tilde\alpha=2\tilde\alpha_i$ for some $\alpha_i\in \pi^*$. (The 
    latter case can occur only if   $\Sigma$ is of type BC$_r$ for some 
    $r$ and $\tilde\alpha/2$ is also a root.)  By 
    the discussion preceding the lemma, there exists $w\in W$ such 
    that $w=\tilde w$ upon restriction to $Q(\Sigma)$.    
     Since  $P(\pi)$ is invariant under the action of $W$, we may assume 
    that $\tilde\alpha$ is either $\tilde\alpha_i$ or 
    $2\tilde\alpha_i$ for some $\alpha_i\in \pi^*$.  
    
   By (1.1), we
can write $$\tilde\alpha_i=(\alpha_{i}+\alpha_{{\rm p}(i)}+\sum_{\alpha_j\in 
\pi_{\Theta}}r_j\alpha_j)/2$$ for some nonnegative integers
$r_j$.     Let
$\omega=\sum_jm_j\omega_j$ be an element of $P(\pi)$.   Note that 
$(\tilde\gamma,\tilde\beta)=(\gamma,\tilde\beta)$ for any weights 
$\gamma$ and $\beta$.
It follows that 
$$(\tilde\omega,\tilde\alpha_i)=(m_ir_i(\alpha_i,\alpha_i)+ 
m_{{\rm p}(i)}r_{{\rm p}(i)}(\alpha_{{\rm p}(i)},\alpha_{{\rm p}(i)})+\sum_{\alpha_j\in 
\pi_{\Theta}}m_jr_j(\alpha_j,\alpha_j))/4.$$  Note that if 
$\alpha_i$ satisfies the conditions of Case 1, then $r_j=0$ for all 
$\alpha_j\in \pi_{\Theta}$ and $\alpha_i=\alpha_{{\rm p}(i)}$.   Hence
$(\tilde\omega,\tilde\alpha_i)=m_ir_i(\alpha_i,\alpha_i)/2\in {\bf 
Z}(\tilde\alpha_i,\tilde\alpha_i)/2.$  Now assume $\alpha_i$ 
satisfies the conditions of Case 2 or Case 3. From the list of  
  rank one symmetric pairs given above, one checks that  $\gg_i$ is not of type G2.  
Hence for  $\alpha_j\in 
\pi_{\Theta}$, $(\alpha_j,\alpha_j)=s(\alpha_i,\alpha_i)$ where 
$s=1,2$, or $1/2$.  Moreover, this last possibility occurs if and 
only if $\gg_i$ is of type $B_r$ as in (3.4)   and $m_j$ equals $2$.   
Hence 
$\sum_{\alpha_j\in 
\pi_{\Theta}}m_jr_j(\alpha_j,\alpha_j)\in {\bf 
Z}(\alpha_i,\alpha_i).$  Since ${\rm p}$ induces a diagram 
automorphism on $\pi$, we have that $(\alpha_{{\rm 
p}(i)},\alpha_{{\rm p}(i)})=(\alpha_i,\alpha_i)$.
Thus $$(\tilde\omega,\tilde\alpha_i)\in {\bf 
Z}(\alpha_i,\alpha_i)/4.\eqno{(3.9)}$$   It follows that 
$(\tilde\omega,\tilde\alpha_i)\in {\bf 
Z}(\tilde\alpha_i,\tilde\alpha_i)/2$ for all 
$\alpha_i$ which satisfies the conditions of Case 2.  Now assume 
that    $\alpha_i$ satisfies the conditions of Case 3.   In 
this case, $(\tilde\alpha_i,\tilde\alpha_i)=(\alpha_i,\alpha_i)/4$.  
Therefore (3.9) implies that 
 $(\tilde\omega,2\tilde\alpha_i)\in {\bf 
Z}(2\tilde\alpha_i,2\tilde\alpha_i)/2.\Box$

\medskip
An immediate consequence of  Lemma 3.3 is that $\widetilde 
{P^+(\pi)}$ is a subset of $P^+(\Sigma)$.  
  The next result will be used to determine which irreducible symmetric pairs satisfy 
$\widetilde {P^+(\pi)}=P^+(\Sigma)$.

\begin{lemma}
    Suppose that each $\alpha_i\in \pi^*$ satisfies one of the 
    following conditions.
   \begin{enumerate}
       \item[(i)] $\Theta(\alpha_i)=-\alpha_i$
       \item[(ii)] $i\neq {\rm p}(i)$
       \item[(iii)]
       $\Theta(\alpha_j)=-\alpha_{{\rm p}(j)}$ for all $j\neq i$
    \end{enumerate}
    Then $\widetilde{P^+(\pi)}=P^+(\Sigma)$.
\end{lemma}

\noindent
{\bf Proof:}  If $\alpha_i$ satisfies (i) or (ii) then Lemma 3.1 
implies that 
 $\tilde\omega_i=\omega_i'$. Suppose 
$\alpha_i$ satisfies the conditions of (iii).  Then $\pi_{\Theta}\cap 
\{\alpha_k|\ (\omega_k,\tilde\alpha_j)\neq 0\}$ is the empty set for $j\neq i$.  Hence 
$(\omega_k,\tilde\alpha_j)=0$ for $\alpha_k\in \pi_{\Theta}$ and 
$j\neq i$.   Let $\lambda_i\in P^+(\pi)$ be chosen as in Lemma 3.3.  It 
follows that $\tilde\lambda_i=\omega_i'$. $\Box$

\medskip

In Theorem 3.6, we find necessary and sufficient conditions for 
$\widetilde{P^+(\pi)}$ to equal $P^+(\Sigma)$.   First, we take a 
look at
a family of symmetric pairs for which this equality fails.
Consider first the  special case when the symmetric pair 
$\gg,\gg^{\theta}$ is of type AII. 
 It follows that 
$\pi_{\Theta}=\{\alpha_{2i+1}|\ 0\leq i\leq t\}$ and the simple 
restricted roots are 
$$\tilde\alpha_{2i}=(\alpha_{2i-1}+2\alpha_{2i}+\alpha_{2i+1})/2\eqno{(3.10)}$$ for 
$1\leq i\leq t$. 
Now assume that the  Lie algebra $\gg$ is of classical type, $\gg$ contains a $\theta$ invariant Lie subalgebra $\gr$ 
such that $\gr,\gr^{\theta}$ is 
of type AII,   and the rank of the restricted root system associated to 
$\gr, \gr^{\theta}$ is $t-1$.  
It follows that   the first $t-1$ simple 
restricted roots    are given by the 
formula in (3.10). We list below the possibilities for 
$\gg,\gg^{\theta}$ under these assumptions.
The formula 
of the last simple restricted 
root is given and   the corresponding rank one symmetric pair is 
identified for the last four types. 
\begin{enumerate}
    \item[(3.11)] $\gg,\gg^{\theta}$ is of type AII and $\Sigma$ is of 
    type A$_t$.  
    \item[(3.12)] $\gg,\gg^{\theta}$ is of type CII(i), $\Sigma$ is 
    of type B$_t$,  
    $\tilde\alpha_{2t}=
    (\alpha_{2t-1}+2\alpha_{2t}+2\alpha_{2t+1}+\cdots+ 
    2\alpha_{n-1}+\alpha_n)/2$, and the rank one symmetric pair 
    corresponding to $\tilde\alpha_{2t}$ is given by (3.7).
    \item[(3.13)] $\gg,\gg^{\theta}$ is of type CII(ii), $\Sigma$ 
    is of type C$_t$,  $\tilde\alpha_n=\alpha_n+\alpha_{n-1}$ and 
    and the rank one symmetric pair 
    corresponding to $\tilde\alpha_{n}$ is given by (3.5).
    \item[(3.14)] $\gg,\gg^{\theta}$ is of type DIII(i), $\Sigma$ is 
    of type C$_t$, $\tilde\alpha_n=\alpha_n$, and the rank one symmetric pair 
    corresponding to $\tilde\alpha_{n}$ is given by (3.1).
    \item[(3.15)] $\gg,\gg^{\theta}$ is of type DIII(ii), $\Sigma$ is 
    of type B$_t$,  $\tilde\alpha_n=(\alpha_n+\alpha_{n-1})/2$,
    and the rank one symmetric pair 
    corresponding to $\tilde\alpha_{n}$ is given by (3.2).
    \end{enumerate}
The next lemma describes the images of the fundamental weights 
inside $P(\pi)$ under $\ \widetilde{\ }$.
We use the following notation just for   the  symmetric pairs listed 
in (3.11)-(3.15) to make the proofs concerning these pairs easier to 
read.
Write $\mu_1,\dots, \mu_t$ 
for the set of positive simple roots in $\Sigma$.  Thus, 
$\{\mu_1,\dots, \mu_t\}=\{\tilde\alpha_i|\ \alpha_i\in \pi^*\}$.    
Similarly, for $1\leq i\leq t$,  set $\eta_i$ equal to the fundamental weight in 
$P(\Sigma)$ corresponding to the simple root $\mu_i$. We further  set $\eta_0=0$
and $\eta_{t+1}=0$.

\begin{lemma} Assume that $\gg,\gg^{\theta}$ satisfies the conditions 
of one of 
(3.11)-(3.15).  Set $t'=t$ if $\gg,\gg^{\theta}$ is of type AII and 
if $\gg,\gg^{\theta}$ is of type 
CII(ii), and 
set $t'=t-1$ otherwise. Then  for all $1\leq i\leq t'$ and $ 0\leq j\leq t'-1$ we have
\begin{enumerate}
    \item[(i)] $\tilde\omega_{2i}=2\eta_i$
    \item[(ii)]$\tilde\omega_{2j+1}=\eta_j+\eta_{j+1}$
    \item[(iii)] If $\gg,\gg^{\theta}$ is not of type CII(ii) then 
    $\tilde\omega_n=\eta_t$.
    \end{enumerate}
\end{lemma}    
    \noindent
    {\bf Proof:}  The first assertion follows from Lemma 3.1 while 
    the second assertion is a straightforward computation using the 
    description of the simple restricted roots in (3.10)-(3.16). The 
    last assertion follows from  Lemma 3.1,
    Lemma 3.2, and 
    the rank one information in the list (3.10)-(3.16).$\Box$
    \medskip

    Assume again that $\gg,\gg^{\theta}$ is of type AII with $t\geq 
    7$. It follows from the previous lemma that $\eta_i\notin 
    \widetilde{P^+(\pi)}$ for $1<i<t$. On the other hand, if 
    $\gg,\gg^{\theta}$ is one of the cases listed in (3.12)-(3.16) and 
    $t\geq 7$, then $\eta_2\notin \widetilde{P^+(\pi)}$. Now suppose that  $\gg,\gg^{\theta}$ is 
of type CII(ii) and $t\geq 4$.  Then one checks that 
$\tilde\omega_n=2\eta_t$.  Furthermore, by Lemma 3.5(ii),
 $\tilde\omega_{n-1}=\eta_{t-1} +\eta_t$.  Thus $\eta_t\notin
\widetilde{P(\pi)}$.  This information is used in 
  the next theorem which determines when $\widetilde{P^+(\pi)}$ 
    equals $P^+(\Sigma)$.

\begin{theorem}  $\widetilde{P^+(\pi)}=P^+(\Sigma)$ if and only if 
$\gg,\gg^{\theta}$ is not of type EIII, EIV, EVII, EIX, or CII(ii), and $\gg$ does 
not contain a $\theta$ invariant
Lie subalgebra $\gr$  of rank greater than or equal to $7$
such that $\gr, \gr^{\theta}$ is of type AII.
\end{theorem}

\noindent
{\bf Proof:} First assume that $\gg,\gg^{\theta}$ is not of type EIII, 
EIV, EVII, EIX, or CII(ii), and $\gg$ does 
not contain a $\theta$ invariant
Lie subalgebra $\gr$  of rank greater than or equal to $7$
such that $\gr, \gr^{\theta}$ is of type AII.   If 
$\gg,\gg^{\theta}$ is not of type EVI, it is straightforward to check that 
each simple 
root  of $\gg$ satisfies the conditions of Lemma 3.4.   On the other 
hand, if $\gg,\gg^{\theta}$ is of type EVI, then 
$\tilde\omega_1=\omega_3'$, $\tilde\omega_2=\omega_3'$, 
$\tilde\omega_5=\omega_5'$ and $\tilde\omega_6=\omega_6'$.   Hence 
$\widetilde{P^+(\pi)}=P^+(\Sigma)$ in this case as well.

Now suppose that $\gg,\gg^{\theta}$ is of type EIV, EVII, or 
EIX.  A straightforward argument  shows that $\omega_1'\notin\tilde{\pi}$.
Hence $\widetilde{P^+(\pi)}$ is a proper subset of $P^+(\Sigma)$ in 
these cases.  Consider the case when $\gg,\gg^{\theta}$ is of type 
EIII.   Then  
$\Theta(\alpha_2)=-\alpha_2-2\alpha_4-\alpha_3-\alpha_5$ and 
$(\Theta(\alpha_2), \alpha_2)=0$.   Hence by Lemma 3.1, 
$\tilde\omega_2=2\omega_2'$.   On the other hand, 
$(\omega_j,\tilde\alpha_1)\neq 0$ for $j\neq 2$.  Thus 
$\omega'_2\notin\{\tilde\omega_1,\dots, \tilde\omega_n\}.$ It follows 
that $\widetilde {P^+(\pi)}$ is not 
equal to $P^+(\Sigma)$ in this case as well.  

The remaining cases follow from the discussion preceding the theorem.
$\Box$

\medskip
Suppose that $P^+(\Sigma)=\widetilde{P^+(\pi)}$.  It follows 
from   Lemma 3.3   that 
$P(\Sigma)=\widetilde{P(\pi)}$.  It turns out that this second equality 
holds for all symmetric pairs.   

    \begin{theorem}  The set $\widetilde{P(\pi)}$ equals $P(\Sigma)$.
	\end{theorem}
	
	\noindent 
	{\bf Proof:} By  Lemma 3.3, we only need to show that 
	$P(\Sigma)$ is a subset of $\widetilde{P(\pi)}$. Suppose that 
	$\widetilde{P^+(\pi)}=P^+(\Sigma)$.  Then 
	the lemma follows from the fact that $\Sigma\subseteq 
	\widetilde{\Delta}$.  Hence, by  Theorem 3.6, we may assume 
	that $\gg,\gg^{\theta}$ falls into one of the following two cases:
       \begin{enumerate}
	  \item[(i)] $\gg,\gg^{\theta}$ is one of the four exceptional 
	symmetric pairs EIII, EIV, EVII, and EIX
       \item[(ii)]  $\gg$ is classical 
	and $\gg$ contains a $\theta$ invariant Lie subalgebra $\gr$  such 
	that $\gr, \gr^{\theta}$ is of type AII  and the rank of $\Sigma$ 
	greater than or equal to $4$.
	\end{enumerate}
	Assume $\gg,\gg^{\theta}$ satisfies the conditions of Case (ii). 
		By Lemma 3.5, $\widetilde{P(\pi)}$ contains $\eta_1$ and 
	$\eta_i+\eta_{i+1}$ for $1\leq i\leq t-2$. Using an inductive 
	argument, we may assume that $\eta_j\in \widetilde{P(\pi)}$ 
	where $j$ is an integer between $1$ and $t-3$.  Then 
	$\eta_{j+1}=(\eta_{j+1}+\eta_j)-\eta_j$ and so $\eta_{j+1}\in 
	\widetilde{P(\pi)}$.    Thus $\eta_1,\dots, \eta_{t-1}$ are all 
	contained in $\widetilde{P(\pi)}$.  Now if $\gg,\gg^{\theta}$ is of type 
	AII or CII(ii), then by Lemma 3.5, $\eta_{t-1}+\eta_{t}$ is also in 
	$\widetilde{P(\pi)}$.   Thus since $\eta_{t-1}\in 
	\widetilde{P(\pi)}$, so is $\eta_t$.  On the other hand if 
	$\Sigma$ is not of type AII or CII(ii), then by Lemma 3.5(iii), 
	$\eta_t\in \widetilde{P(\pi)}$.  This completes the proof in 
	Case (ii).
	
	Consider the cases when $\gg,\gg^{\theta}$ is of type EIV, 
	EVII, or EIX.   
    A straightforward computation shows that   
    $\omega_1'=\tilde\omega_5-\tilde\omega_6$ and 
    $\omega_6'=\tilde\omega_3-\tilde\omega_1$.  Furthermore, by  
    Lemma 3.1, we have $\tilde\omega_i=\omega_i'$ for $i>7$.   Thus 
    $\widetilde{P(\pi)}$ contains all the fundamental restricted roots 
    which imply the lemma for these cases.  Now assume $\gg,\gg^{\theta}$ 
    is of type EIII.  By  Lemma 3.1, we have $\tilde\omega_1=\omega_1'$.
    One further checks that $\omega_2'=\tilde\omega_3-\tilde\omega_1$ 
    which completes the proof of the theorem. $\Box$

\section{Quantum Radial Components}

Throughout this section, we assume that  $B$ be a fixed coideal subalgebra in ${\cal B}$.  In this section, we study the   map ${\cal X}$ on the 
algebra $\check U^B$ which computes the quantum radial components.   
The map  ${\cal X}$ is
  defined in [L5, Section 3] for standard analogs associated to symmetric pairs 
with reduced root systems.   After extending the definition   to  the general 
case,  we 
 explore connections between ${\cal X}$ and the 
Harish-Chandra projection map ${\cal P}_B$.   In particular, we show that ${\cal 
X}(\check U^B)$ and ${\cal P}_B(\check U^B)$ are isomorphic as 
algebras. As in [L5, Theorem 3.4], one 
has that the   quantum radial components are invariant under the 
action of the restricted Weyl group $W_{\Theta}$.   Using this fact,
we determine the  possible  top degree terms   of   
 elements in ${\cal P}_B(\check U^B) $   with respect to a certain 
 filtration of $\check U$.  This is a crucial step   towards finding the 
 image of $\check U^B$ under the quantum Harish-Chandra map ${\cal 
 P}_B$.

Before defining the function ${\cal X}$, we review some notions from 
[L5, Section 3] 
which are necessary to define  the codomain of this function.
  Both group rings  
${\cal C}[Q(\Sigma)]$ and   ${\cal C}[{\cal A}]$
embed in the right endomorphism ring     ${\rm End}_r{\cal 
C}[P(\Sigma)]$  of ${\cal C}[P(\Sigma)]$ as follows.  The ring ${\cal C}[Q(\Sigma)]$ acts on
$ {\cal 
C}[P(\Sigma)]$   by right multiplication. The action of
${\cal C}[{\cal A}]$ on $ {\cal 
C}[P(\Sigma)]$ is given by 
$$z^{\lambda}*\tau(\mu)=q^{(\lambda,\mu)}z^{\lambda}$$
for all $\lambda\in P(\Sigma)$ and $\tau(\mu)\in {\cal A}$.
Denote the  subring of ${\rm End}_r{\cal C}[P(\Sigma)]$ 
generated by ${\cal C}[Q(\Sigma)]$ and ${\cal A}$  by 
${\cal C}[Q(\Sigma)]{\cal A}$.  One can localize 
by the nonzero elements of ${\cal C}[Q(\Sigma)]$ to form a new ring 
${\cal C}(Q(\Sigma)){\cal A}$ generated by ${\cal A}$ 
and the quotient field ${\cal C}(Q(\Sigma))$ of ${\cal C}[Q(\Sigma)]$.

Now ${\cal C}[{\cal A}]$ is a left ${\cal C}[Q(\Sigma)]{\cal 
A}$ 
module where elements of ${\cal A}$ act via left 
multiplication and
$$z^{\lambda}\cdot \tau(\mu)=q^{(\lambda,\mu)}\tau(\mu) $$ for all 
$\lambda\in Q(\Sigma)$ and $\tau(\mu)\in {\cal A}$. Note that $g\cdot \tau(\mu)$ is just a scalar 
multiple of $\tau(\mu)$ for $\tau(\mu)\in {\cal A}$ and $g\in {\cal 
C}[Q(\Sigma)]$. Given an element $f\in {\cal 
C}[Q(\Sigma)]{\cal A}$ and $g\in {\cal C}[Q(\Sigma)]$, set
$$(fg^{-1})\cdot \tau(\mu)=(f\cdot \tau(\mu))((g\cdot 
\tau(\mu))\tau(\mu)^{-1})^{-1}$$ for all $\tau(\mu)\in {\cal A}$ such 
that 
$(g\cdot 
\tau(\mu))$ is nonzero.

The restricted Weyl group  $W_{\Theta}$ acts on   ${\cal 
C}[Q(\Sigma)]{\cal A}$ where $w\tau(\beta)=\tau(w\beta)$ and 
$w(z^{\mu})=z^{w\mu}$ for all $w\in W_{\Theta}$, $\tau(\beta)\in 
{\cal A}$ and $\mu\in Q(\Sigma)$. This action of $W_{\Theta}$ on   ${\cal 
C}[Q(\Sigma)]{\cal A}$  extends to an action of $W_{\Theta}$ on 
${\cal C}(Q(\Sigma)) {\cal A}$.

Recall that $B$ is a fixed coideal subalgebra of ${\cal B}$.  Assume 
that $B'$ is another subalgebra of ${\cal B}$.      Given $\lambda\in P^+(2\Sigma)$, let $g_{\lambda}$ denote 
the zonal spherical function at 
$\lambda$.  In particular, $g_{\lambda}$ is a basis vector for the one 
dimensional  space of 
right $B$ and left $B'$ invariants in $L(\lambda)\otimes 
L(\lambda)^*$. (For more information about zonal spherical functions, 
see [L4] and [L5].)
Since elements of $L(\lambda)\otimes L(\lambda)^*$ can be viewed as 
functions on $U$, restriction to $\check U^0$ defines a map from   $L(\lambda)\otimes L(\lambda)^*$
to ${\cal C}[P(\pi)]$.   Write $\varphi_{\lambda}$ for the image of 
$g_{\lambda}$ in ${\cal C}[P(\pi)]$ under this restriction map.   
By [L4, Theorem 4.2],  $\varphi_{\lambda}$ is  an element of  ${\cal 
C}[P(2\Sigma)]$.   We further assume that $B'$ has been chosen so 
that $\varphi_{\lambda}$ is $W_{\Theta}$ invariant ([L4, Theorem 5.3]). The 
observant reader might notice that $B'$ is not mentioned again in the 
results and discussion below.   In addition to being used to specify 
the zonal spherical functions $g_{\lambda}$, $B'$ is also used to define the map 
${\cal X}$ ([L5, Section 3]).  So officially, both $g_{\lambda}$ and ${\cal X}$ 
depend on the choice of $B'$ once $B$ is fixed.  However,  
the information we get in terms of the Harish-Chandra map ${\cal P}_B$ 
does not depend at all on $B'$.

The following theorem is a generalization of [L5, Theorem 3.6].

\begin{theorem}
    There is an algebra homomorphism  $${\cal X}:\check U^B\rightarrow
   ({\cal C}(Q(\Sigma)){\cal A})^{W_{\Theta}}$$ such that 
    $$\varphi_{\lambda}*{\cal X}(c)=z^{\lambda}({\cal P}_B(c))\varphi_{\lambda}
    \quad{\rm  and  }\quad
    g_{\lambda}(c\tau(\beta))=(\varphi_{\lambda}*{\cal 
    X}(c))(\tau(\beta))$$ for all $c\in \check U^B$, $\lambda\in 
    P^+(2\Sigma)$, and $\tau(\beta)\in {\cal A}$.
    \end{theorem}
    
\noindent
{\bf Proof:} The results and constructions of [L5, Section 3] 
generalize to all coideal subalgebras in ${\cal B}$ using    Lemma 2.1 and Theorem 
 2.2 of this paper
instead of [L5, Lemma 2.1] and [L5, Theorem 2.2].  Thus the  existence 
of  a map 
$${\cal X}:\check U^B\rightarrow
   ({\cal C}(Q(\Sigma))\check{\cal A})^{W_{\Theta}}$$ which satisfies 
   the conclusions of the theorem follows from the same arguments as in   
   the proof of  
  [L5, Theorem 3.6].   Lemma 2.5 and the last assertion 
   of [L5, Theorem 3.2] ensure that the image of $ \check U^B$ 
   under ${\cal X}$ is a 
   subset of ${\cal C}(Q(\Sigma)){\cal A}$.
   $\Box$
   
   \medskip
   
   In [L5, Section 5], a filtration ${\cal F}$ on $U$ is introduced
so that $B\cong {\rm  gr}_{\cal F}B$ whenever $B$ is a standard analog and 
the restricted root system is reduced ([L5, Lemma 5.1]).  The filtration 
is defined using a degree function defined by 
\begin{enumerate}
    \item[(4.1)]$\deg x_i=\deg y_it_i=0$ for all $1\leq i\leq n$
    \item[(4.2)] $\deg t_i^{-1}=1$ for all $\alpha_i\in \pi\setminus \pi_{\Theta}$
    \item[(4.3)]  $\deg t=0$ for all $t\in T_{\Theta}$.  
 \end{enumerate}
  This filtration can be further extended to $\check U$ as follows.  
 Given $\mu\in P(\pi)$ such that $\tilde\mu=\sum_im_i\tilde\alpha_i$, 
 we set  
  $\deg(\tau(\mu))=\sum_im_i\tilde\alpha_i$.  As explained in the 
  proof of [L5, Lemma 5.1], $\deg \tilde\theta(y_i)t_i=0$ for all 
  $\alpha_i\in \pi\setminus\pi_{\Theta}$.   Since $\deg t_i=-1$ for 
  $\alpha_i\in \pi\setminus\pi_{\Theta}$, it follows that $\deg B_i=\deg 
  (y_it_i+d_i\tilde\theta(y_i)t_i+s_it_i=0$ for all $\alpha_i\in 
  \pi\setminus\pi_{\Theta}$.  The arguments of [L5, Lemma 5.1] can thus 
  be used to show
  that  $B\cong {\rm gr}_{\cal F}B$ 
  for all $B\in {\cal B}$.   The rest of [L5, Section 5] 
  extends to the nonreduced restricted root system and nonstandard 
  analog cases verbatim.   The proof of   [L5, 
  Theorem 5.8 and Corollary 
6.8] shows that there exists a  formal power series $p$ in the $z^{-\tilde\mu}$,
$\tilde\mu\in Q^+(\Sigma)$, such that 
    $${\rm gr}_{\cal F}({\cal X}(u))=p^{-1}{\rm gr}_{\cal 
F}({\cal P}_B(u))p  \eqno{(4.4)}$$ for all $u\in \check U^B$.  Here we are viewing the ring ${\cal 
C}(Q(\Sigma)){\cal A}$ as embedded in the larger ring generated by    formal 
Laurent series in the $z^{-\tilde\eta}$, $\tilde\eta\in Q(\Sigma)$,
and elements of ${\cal A}$ (see [L5, Section 3, discussion preceding 
(3.9)]).   

It should be noted that an   explicit formula 
for $p$  is provided by  [L5, Lemma 6.6] for standard quantum analogs 
with reduced restricted 
root systems.   However, this formula   is not needed in the 
proof of the first assertion of [L5, Corollary 6.8] which is just (4.4) 
above.  Indeed, with the exception of [L5, Lemma 6.6], the results of 
Section 6 also carry over verbatim to all irreducible symmetric pairs 
and analogs in ${\cal B}$. The results of [L5] and Theorem 4.1 yield 
the following isomorphism of 
 rings.

 \begin{corollary}
     The map ${\cal X}(u)\rightarrow {\cal P}_B(u)$ defines an algebra 
     isomorphism from ${\cal X}(\check U^B)$ onto ${\cal P}_B(\check U^B)$.
     Moreover, the kernel of both ${\cal X}$ and ${\cal P}_B$ upon 
     restriction to $\check U^B$ is 
     equal to $\check U^B\cap (B\check T_{\Theta})_+\check U$. 
      \end{corollary}
   
      \noindent
      {\bf Proof:} Suppose that $u\in \check U^B$.
        Theorem 4.1 ensures that ${\cal P}_B(u)=0$ if and 
 only if ${\cal X}(u)=0$. This proves the  first assertion.  The 
 second assertion follows from  [L5, Corollary 6.8] and the comments 
 at the end of [L5, Section 6].   
 $\Box$

\medskip
  
Note that the filtration ${\cal F}$ restricts to a filtration on 
${\cal C}[{\cal A}]$.  Moreover, the degree function which defines 
${\cal F}$ actually gives ${\cal C}[{\cal A}]$ the 
structure of a graded algebra.   Given $X\in {\cal
C}[{\cal A}]$, let $\tip(X)$ be the element of ${\cal
C}[{\cal A}]$ homogeneous of degree $\deg(X)$ such that $X-\tip(X)$
has degree strictly less than $\deg(X)$.   We can further extend 
${\cal F}$ to a  filtration on ${\cal C}(Q(\Sigma)){\cal A}$ by insisting that 
nonzero elements of ${\cal C} (Q(\Sigma))$ have degree $0$.  Moreover, 
  the definition of $\tip$ can be extended to elements of ${\cal C}(Q(\Sigma)){\cal A}$.

The next lemma uses the close connection between ${\cal X}$ and 
${\cal P}_B$ to gain further information about the image of $\check 
U^B$ under ${\cal P}_B$.

\begin{lemma}The set $\{\tip(X)|X\in{\cal P}_B(\check U^B)\}$ is a 
subset of the ${\cal C}$ span of $\{\tau(-2\eta)|\eta\in P^+(\Sigma)\}.$
\end{lemma}

\noindent
{\bf Proof:}   Let $u\in {\cal P}_B$. By Theorem 4.1, we can write 
$$\tip({\cal X}(u))=
\sum_{\beta\in S(u)}\tau(2\beta)f_{\beta}$$ for some finite subset 
$S(u)$ of $P(\Sigma)$ and nonzero elements $f_{\beta}\in {\cal C}(Q(\Sigma))$. 
Theorem 4.1 ensures that ${\cal X}(u)$ is invariant under 
the action of $W_{\Theta}$. Now each $W_{\Theta}$ orbit in 
$P(2\Sigma)$ contains an antidominant weight.  Moreover, this 
antidominant weight is smaller than every other element in its $W_{\Theta}$ orbit with respect to the 
standard partial order on $\gh^*$.       Applying $\tau$  yields 
the corresponding result for $W_{\Theta}$ orbits in ${\cal A}$.  It 
follows that      $-\beta\in P^+(\Sigma)$ for each $\beta\in S(u)$. The lemma now follows from the 
expression for ${\rm gr}_{\cal F}({\cal X}(u))$ given in (4.4). $\Box$

\section {The ordinary Harish-Chandra map}

  In this section, we begin the study of the center of 
$\check U$ under the map ${\cal P}_B$, for $B\in {\cal B}$.   The main tool here is the 
  the ordinary quantum Harish-Chandra map, ${\cal P}$. We see below 
  that it is easy to compute ${\cal P}_B(z)$ once ${\cal P}(z)$ is 
  known.   Thus we analyze ${\cal P}_B(Z(\check U))$ by 
  exploiting detailed information concerning the image 
  of a basis of central elements under the ordinary Harish-Chandra map ${\cal P}$. 
  This information is used to   show that the image of the center $Z(\check 
U)$ under the map ${\cal P}_B$ is invariant under a dotted action of 
the restricted Weyl group $W_{\Theta}$.    As we will see in Section 
6, this is enough to show that 
  ${\cal P}_B(\check U^B)$ is 
  $W_{\Theta}.$ invariant for almost all symmetric pairs $\gg,\gg^{\theta}$.
 The material in this section and the fact that     ${\cal P}_B(\check U^B)$ is 
  ``large'' inside of ${\cal P}_B(\check U^B)$ is used to extend this invariance 
  result  
  to the remaining cases in [L6].

We recall here the definition of the 
ordinary quantum  Harish-Chandra projection. Using the triangular decomposition of 
 $\check U$, $\check U$ admits the following direct sum decomposition.
$$\check U=\check U^0\oplus(G^-_+\check U+\check UU^+_+).\eqno{(5.1)}$$ Let 
${\cal P}$ denote the projection of $\check U$ onto $\check U^0$ using 
this decomposition.  Let $B\in {\cal B}$.  Note that  ${\cal P}_{B}\circ 
{\cal P}$ maps $\check U$ onto ${\cal C}[\check{\cal A}]$.   In 
general, ${\cal P}_{B}\circ {\cal P}$ does not agree with ${\cal 
P}_{B}$.   However, the discussion preceding     [L5, Theorem 3.6] 
shows that the restriction of ${\cal P}_{B}$ to $Z(\check U)$ 
is the same function as the restriction of   
${\cal P}_B\circ{\cal P}$ to $Z(\check U)$. Now 
${\cal P}_B\circ{\cal P}$ is the same as the composition of ${\cal P}$ 
followed by the projection of $\check U^0$ onto ${\cal C}[\check {\cal 
A}]$ using (2.3).  In particular, the restriction of ${\cal P}_B$ to 
$Z(\check U)$ is {\it independent} of the choice of $B\in {\cal B}$. Moreover, 
$z^{\lambda}({\cal P}(z))=z^{\lambda}({\cal P}_B(z))$ for all $z\in 
Z(\check U)$ and $\lambda\in P^+(2\Sigma)$.

  A description of a nice basis for the  center of $Z(\check U)$ and the image 
  of this basis under ${\cal P}$  is given in [Jo, Chapter 7] and [JL].  It should 
  be noted that these references
  use the locally finite part of $\check U$ with respect to the left 
  adjoint action.    Since we are viewing the center as a subalgebra of the 
  $(\adr B)$ invariants of $\check U$, it is necessary to translate
  these results  to the setting of the 
   right 
  adjoint action.   
In particular, for each   $\mu\in P^+(\pi)$ there exists  a unique central element 
$z_{2\mu}$ in 
$\tau(2\mu)+(\adr U_+)\tau(2\mu)$. 
Moreover, the set 
$\{z_{2\mu}|\mu\in P^+(\pi)\}$ forms a basis for $Z(\check U)$.
Let $\rho$ denote the half sum of the positive roots in $\Delta$.  Given $\lambda\in P^+(\pi)$, set $$\hat 
\tau(\lambda)=\sum_{w\in W} 
\tau(w\lambda)q^{(\rho,w\lambda)}.\eqno{(5.2)}$$  The 
next lemma, which  describes the image of each of these basis elements 
under ${\cal P}$, is a version of [Jo, Lemma 7.1.19] with respect to  the right adjoint action.

\begin{lemma} For all $\mu\in P^+(\pi)$, $${\cal P}(z_{2\mu})=a_{2\mu}^{-1}\sum_{\nu\in 
P^+(\pi)}\hat \tau(2\nu)\dim L(\mu)_{\nu}$$ where 
$$a_{2\mu}= \sum_{\nu\in 
P^+(\pi)}(\sum_{w\in W} q^{(\rho,2w\nu)})\dim L(\mu)_{\nu}. $$

    \end{lemma}
    
    \noindent
    {\bf Proof:}  Let $\iota$ be the ${\bf C}$ algebra involutive 
    antiautomorphism of $\check U$ 
    defined by $\iota(x_i)=y_i$ and $\iota(y_i)=x_i$ for all $1\leq 
    i\leq n$, $\iota(q)=q^{-1}$, and $\iota(t)=t^{-1}$ for all $t\in 
    \check U$.  
  It is straightforward to show that $$\iota((\ad 
    a)b)=-(\adr \iota(a))\iota(b)\eqno{(5.3)}$$ for all $a\in \check U$ and $b\in 
    \{x_i,y_i|\ 1\leq i\leq n\}\cup \check T.$  Hence (5.3) holds for 
    all $a\in \check U$ and $b\in U$.   It follows that $$\iota((\ad 
    U_+)\tau(-2\mu))=(\adr U_+)\tau(2\mu)$$ for all $\mu\in P^+(\pi)$. 
    Thus $\iota(z_{2\mu})$ is the unique central element contained 
    in $\tau(-2\mu)+(\ad U_+)\tau(-2\mu)$.  Since $\iota$ preserves 
    both $\check U^0$ 
    and $G^-_+\check U+\check UU^+_+$, it follows that   $\iota$ 
    commutes with the   map ${\cal P}$.   In particular ${\cal 
    P}(z_{2\mu})=\iota({\cal P}(\iota(z_{2\mu})))$.  Thus by [Jo, Lemma 7.1.19], 
    we have
    $${\cal P}(z_{2\mu})=a_{2\mu}^{-1}\sum_{\nu\in 
P^+(\pi)}\hat\tau(2 \nu)\dim L(\mu)_{\nu}$$ 
where $a_{2\mu}$ is a nonzero scalar.

Recall the information about zonal spherical functions discussed 
right before Theorem 4.1.
The spherical function at $0$ is just the basis vector $g_0$ for 
the one-dimensional module $L(0) \otimes L(0)^*$. Rescaling if 
necessary, we may further assume that $g_0(1)=1$.  Since $\check U_+$ 
annihilates $L(0)$, we have
$g_0(\check U_+)=0$. Now $\tau(2\mu)\in 1+\check U_+$ and $(\adr 
U_+)\tau(2\mu)\in \check U_+$.  It follows that $g_0(z_{2\mu})=g_0(1).$  
On the other hand,  $z^0({\cal 
P}(z_{2\mu}))g_0(1)=z^0({\cal P}_B(z_{2\mu}))=g_0(z_{2\mu}).$  Hence $z^0({\cal 
P}(z_{2\mu}))=1 $ and so $$a_{2\mu}=z^0(\sum_{\nu\in 
P^+(\pi)}\hat \tau(2 \nu)\dim 
L(\mu)_{\nu}).\Box$$

\medskip
Fix $\mu\in P^+(\pi)$. Let $w_0$ denote the 
longest element of $W$. Set $a=a_{-2w_0\mu}$ where $a_{-2w_0\mu}$
is the scalar defined  in the previous lemma for the central element 
$z_{-2w_0\mu}$.  
Using Lemma 5.1, we can write $a{\cal P}(z_{-2w_0\mu})$ as a sum  
$$b_{-\mu}\tau(-2\mu)+\sum_{\gamma>-\mu}b_{\gamma}\tau(2\gamma)$$  
where the $b_{\gamma}$ are scalars. Since $\dim 
L(-w_0\mu)_{-\mu}=1$, it follows that   $b_{-\mu}$ is just a power of 
$q$ times the order of the stabilizer of $-\mu$ in $W$.  More generally, $b_{\gamma}\in {\bf N}[q,q^{-1}]$ 
for each  $\gamma$.
 Now ${\cal 
P}_B(z_{-2w_0\mu})={\cal P}_B\circ {\cal P}(z_{-2w_0\mu})$. 
Hence $$a{\cal P}_B(z_{-2w_0\mu})\in 
(\sum_{\tilde\gamma=\widetilde{\mu}}b_{-\gamma})\tau(-2\tilde\mu)+
\sum_{\tilde\gamma>-\tilde{\mu}}{\cal C}\tau(2\tilde\gamma).$$   It 
follows that $\sum_{\tilde\gamma=-\tilde{\mu}}b_{-\gamma}$ is 
nonzero.  Thus $\tip({\cal P}_B(z_{-2w_0\mu})$ is a nonzero multiple of 
$\tau(-2\tilde{\mu})$.  Hence
$${\rm 
span}\{\tau(-2\tilde\mu)|\ \mu\in P^+(\pi)\}\subseteq\{\tip({\cal P}_B(z))|\ z\in Z(\check U)\}.\eqno{(5.4)}$$

  The next lemma provides a sufficient, but not 
necessary, condition for 
${\cal P}_{B}(Z(\check U))$ to equal  ${\cal P}_{B}(\check U^B)$.
\medskip
 \begin{lemma}  If $\widetilde{P^+(\pi)}=P^+(\Sigma)$ then ${\cal P}_B(Z(\check U))=
     {\cal P}_B(\check U^B)$.
\end{lemma}

\noindent
{\bf Proof:} By Lemma 4.3, the set $\{\tip({\cal P}_B(u))|u\in \check 
U^B\}$ is a subset of the span of the set $\{\tau(-2\tilde\mu)|\tilde\mu\in 
P^+(\Sigma)\}$. The lemma now follow  from (5.4).   $\Box$

\medskip

  Recall that there is a dotted action 
of $W$ on 
$\check U^0$ defined by
$$w.\tau(\mu)q^{(\rho,\mu)}=\tau(w\mu)q^{(\rho, w\mu)}$$ for all
$\tau(\mu)$ in $\check T$ and $w\in W$. 
We use the same formula and notation to   define a dotted action of
$W_{\Theta}$ on
${\cal C}[{\cal A}]$ where now we assume that $\tau(\mu)\in {\cal A}$ 
and $w\in W_{\Theta}$.  By [Jo, 7.1.17 and 7.1.25], 
the image of $Z(\check U)$ under the ordinary Harish-Chandra map 
${\cal P}$ is invariant under the dotted action of $W$.

Let $W'$ be the Weyl group generated by the 
reflections associated to the simple roots in $\pi_{\Theta}$.
 
\begin{lemma}  Suppose that $f\in Z(\check U)$.   Then 
$P_B(f)$ is invariant under the dotted action of the restricted Weyl 
group ${W_{\Theta}}$.
\end{lemma}

\noindent
{\bf Proof:}  Suppose that $\tilde w\in W_{\Theta}$.  Let    $w$ be   an 
element of $W$ such that $\tilde w$ and $w$ agree on $Q(\Sigma)$ (see 
the discussion preceding Lemma 3.3).  In particular, we can view  $w$ as an automorphism of 
$Q(\Sigma)$.  Note that   $\sum_{\tilde\alpha\in 
\Sigma}{\bf Q}\tilde\alpha$ is orthogonal to 
$\sum_{\alpha\in \pi}{\bf Q}(\alpha+\Theta(\alpha))$  with respect to 
the Cartan inner product.   It follows that $w$ restricts to an 
automorphism of  $\sum_{\alpha\in \pi}{\bf 
Q}(\alpha+\Theta(\alpha))$. Note that the intersection of  
 $\sum_{\alpha\in \pi}{\bf 
Q}(\alpha+\Theta(\alpha))$ with $\Delta$  is just the root system 
corresponding to $\pi_{\Theta}$.  In particular, $w$ restricts to an 
automorphism of the root system generated by 
$\pi_{\Theta}$. 
It follows that there exists $w'\in W'$ 
such that $ww'$ restricts to a diagram automorphism of 
$\pi_{\Theta}$.  Now  $(\alpha_i,\mu)=0$ for all $\alpha_i\in 
\pi_{\Theta}$ and $\mu\in P(\Sigma)$.   Thus elements of $W'$ restrict 
to the identity on $P(\Sigma)$.   Hence we  may replace $w$ by 
$ww'$.   In particular, we may assume that $w$ acts the same as $\tilde w$ on $Q(\Sigma)$ 
and $w$ restricts to a  diagram automorphism of $\pi_{\Theta}$.

Now suppose $\mu\in P(2\pi)$.    Note that 
$\mu-\tilde\mu=(\mu+\Theta(\mu))/2$.  Note further that $(\rho, 
\alpha)=(\rho,w\alpha)$ for all $\alpha\in \pi_{\Theta}$ since $w$ is 
a diagram automorphism on $\pi_{\Theta}$.
Hence $(\rho, \mu-\tilde\mu)=0$.   It follows that
$$\eqalign{
w.q^{(\rho,\tilde\mu)}\tau(\mu)&=q^{(\rho,w\mu)}q^{(\rho,-\mu
+\tilde\mu)}\tau(w\mu)\cr
&=q^{(\rho,w\tilde\mu)}\tau(w\tilde\mu)
\tau(w(\mu+\Theta(\mu)))\cr&=q^{(\rho,w\tilde\mu)}\tau(w\tilde\mu)+q^{(\rho,w\tilde\mu)}
\tau(w\tilde\mu)
(\tau(w(\mu+\Theta(\mu)))-1).\cr}$$
Thus ${\cal P}_B(  {w.\tau(\mu)})=\tilde w.\tau(\tilde\mu)$. It 
should be noted here that $w.\tau(\mu)$ refers to the dotted action 
of the big Weyl group $W$ while $\tilde w.\tau(\tilde\mu)$ refers to 
the dotted action of the restricted Weyl group $W_{\Theta}$.  Since 
$w.f=f$, it follows that $\tilde w.P_B(f)=P_B(f)$. $\Box$

\medskip
Given $\tilde\mu\in P^+(\Sigma)$, set $$\hat 
m(2\tilde\mu)=\sum_{\gamma\in 
W_{\Theta}\tilde\mu}q^{(\rho,2\gamma)}\tau(2\gamma).$$   Note that 
$\hat m(2\tilde\mu)$ is an element of ${\cal C}[{\cal 
A}]^{W_{\Theta}.}$. Moreover, the set $\{\hat 
m(2\tilde\mu)|\tilde\mu\in P^+(\Sigma)\}$ forms a basis for ${\cal 
C}[{\cal A}]^{W_{\Theta}.}$. Now $\tip(\hat 
m(2\tilde\mu))=q^{(\rho,2\tilde\mu')}\tau(2\tilde\mu')$ where 
$\tilde\mu'$ is the image of $\tilde\mu$ under the longest element in 
$W_{\Theta}$.   Hence $$\{\tip(u)|\ u\in  {\cal C}[
{\cal A}]^{W_{\Theta}.}\}={\rm 
span}\{\tau(-2\tilde\mu)|\ \mu\in P^+(\Sigma)\}.\eqno{(5.5)}$$
The next theorem is part of a series of results that describe the image 
of $\check U^B$ under the map ${\cal P}_B$.

\begin{theorem} Suppose that $\gg,\gg^{\theta}$ is not of type EIII, EIV, EVII, 
EIX, or CII(ii) and $\gg$ does not contain a $\theta$ invariant Lie subalgebra 
$\gr$ of rank greater than or equal to $7$ such that 
$\gr,\gr^{\theta}$ is of type AII.  Then ${\cal P}_B(Z(\check 
U))={\cal P}_B(\check U^B)$.  Moreover, ${\cal P}_B(\check U^B)={\cal 
C}[{\cal A}]^{W_{\Theta}.}$.
    \end{theorem}
    
    \noindent
    {\bf Proof:}  The first assertion follows from Theorem 3.6 and  
    Lemma 5.2. Furthermore, Theorem 3.6, Lemma 4.3,  and (5.4) implies that  
      $$\{\tip({\cal 
    P}_B(z))|\ z\in Z(\check U)\}={\rm 
span}\{\tau(-2\tilde\mu)|\ \tilde\mu\in P^+(\Sigma)\}.$$ Hence by the 
discussion preceding the theorem, we have 
$$\{\tip({\cal 
    P}_B(z))|\ z\in Z(\check U)\}=\{\tip(u)|\ u\in  {\cal C}[
{\cal A}]^{W_{\Theta}.}\}.$$    Now Theorem 2.6 ensures 
that ${\cal P}_B(Z(\check U))$ is a subring of ${\cal C}[{\cal A}]$.
The theorem now follows from the fact,
guaranteed by Lemma 5.3, that    
${\cal P}_B(Z(\check U))$ is $W_{\Theta}.$ invariant.  $\Box$

\medskip
The proof of Theorem 5.4 uses the fact that $\widetilde 
{P^+(\pi)}=P^+(\Sigma)$ for the symmetric pairs under consideration.   
We see in the next   section    that this is not a necessary condition for 
the conclusions of Theorem 5.4 to hold.   On the other hand, we show below 
that for symmetric pairs $\gg,\gg^{\theta}$ 
  of type EIII, EIV, EVII, and  EIX, we have ${\cal P}_B(Z(\check U))$ 
is a proper subset of ${\cal C}[{\cal A}]^{W_{\Theta}.}$.  In 
particular, the above theorem is not true for these symmetric pairs.  
Despite this failure, we show that the last assertion of the above theorem does hold in 
the remaining four types in [L6, Theorem 4.1].

  Let $\check U_{{\bf C}(q)}$ denote the ${\bf C}(q)$ 
subalgebra of $\check U$ generated by $x_i,y_i$ for $1\leq i\leq 
n$ and $\check T$.  Now the center $Z(\check U)$ is isomorphic to a  polynomial ring 
in $n$ variables.  Moreover, we may choose the generators of 
$Z(\check U)$ to be elements of the smaller algebra $\check U_{{\bf 
C}(q)}$ as in [Jo, Section 7].  In particular, we may restrict our attention to the ${\bf 
C}(q)$ algebra $\check U_{{\bf C}(q)}$. 

The next result is proved using specialization techniques and 
classical results.  First we recall basic notions concerning 
specialization.  Set $A$ equal to the localization ${\bf 
C}[q]_{(q-1)}$ of ${\bf C}[q]$ at 
the maximal ideal $(q-1)$.  Let $\hat U$ denote the $A$ 
subalgebra of $\check U_{{\bf C}(q)}$ generated by $x_i,y_i$ for 
$1\leq i\leq n$, $\check T$, and 
$(t-1)/(q-1)$ for all $t\in \check T$. We have an isomorphism 
$$\hat U\otimes_{A}{\bf C}\rightarrow U(\gg)\eqno{(5.6)}$$
(see for example [L1, Section 2]). Given a subalgebra $S$ of $\check 
U$, we say that $S$ specializes to the subalgebra $\bar S$ of 
$U(\gg)$ provided that the image of $S\cap \hat U$ in  $U(\gg)$ is $\bar S$. 

To make the exposition easier, we assume that $\gg$ is simply laced 
and so all roots have the same length (which we assume is $2$).  Ultimately, our focus will be 
on the four exceptional symmetric pairs where the underlying Lie 
algebra $\gg$ is simply laced of type E$_8$.   So this additional 
assumption does not affect the next theorem.

Let $h_1,\dots, h_n$ be a basis for the Cartan subalgebra $\gh$ of $\gg$ 
where $h_i$ can be identified with the coroot of $\alpha_i$ via the 
Killing form.  We have that $(t_i-1)/(q-1)$ is sent to $h_i$ for all $1\leq 
i\leq n$ under the isomorphism (5.6).   More generally, if 
$\beta=\sum_im_i\alpha_i$, then $(\tau(\beta)-1)/(q-1)$ specializes 
to $\sum_im_ih_i$.  Note further that $t$ specializes to $1$ for 
each $t\in \check T$.  Set $T_2=\{\tau(2\mu)|\ \mu\in P(\pi)\}$.  Note 
that $(\tau(2\alpha_i)-1)(q-1)^{-1}$ specializes to $2h_i$.
It follows that the specialization of ${\cal C}[T_2]$ at $q=1$ is just 
the enveloping algebra $U(\gh)$ of the Cartan 
subalgebra $ \gh$.

Let $\ga$ be the subspace of $\gh$ spanned by the set 
$\{h-\theta(h)|h\in \gh\}$.      Note that $\tilde h=(h-\theta(h))/2$ for 
all $h\in \gh$ defines a projection of $\gh$ onto $\ga$.  Moreover, 
this projection extends in the obvious way to an algebra homomorphism  of 
$U(\gh)$ onto $U(\ga)$.   Now the image of ${\cal C}[{\cal A}]$ under 
specialization is just $U(\ga)$.  Recall that the restriction of 
${\cal P}_B$ to $\check U^0$ sends $\tau(\mu)$ to $\tau(\tilde\mu)$ 
for all $\mu\in P(\pi)$.  Also, Lemma 3.3 
ensures that the image of 
${\cal C}[T_2]$ under ${\cal P}_B$ is a subalgebra of ${\cal 
C}[{\cal A}]$. In particular, we have a commutative 
diagram:
$$\matrix{{\cal C}[T_2]&\rightarrow &U(\gh)\cr
             \ &\ &\ \cr
	     \downarrow&\ &\downarrow\cr
	     \ &\ &\ \cr
	      {\cal C}[{\cal A}]&\rightarrow &U(\ga)}\eqno{(5.7)}$$ 
where the top and bottom maps are specialization and the downward maps 
are algebra homomorphisms defined using $\ \tilde{\ }\ $. Moreover, 
all
maps  in this 
diagram are  surjections.

Let $\alpha$ be a simple root in $\pi$ and let 
$s_{ \alpha}$ denote the corresponding reflection in $W$.
For each $i$ such that $1\leq i\leq n$, we have 
$$\eqalign{s_{ \alpha}\cdot 
{{(t_i-1)}\over{(q-1)}}&={{(q^{(\rho,s_{ \alpha}\alpha_i-\alpha_i)}t_i-1)}\over{ (q-1)}}\cr
&={{(q^{(\rho,s_{ \alpha}\alpha_i-\alpha_i)}-1)}\over{(q-1)}}t_i+{{(t_i-1)}\over{(q-1)}}.\cr}$$
Since $\gg$ is simply laced it follows that 
$s_{\alpha}\alpha_i=\alpha_i-(\alpha_i,\alpha)\alpha$.   
Furthermore, $(\rho, \alpha)=1$.  It follows that $s_{ \alpha}\cdot (t_i-1)(q-1)^{-1}$ specializes to 
$h_i-(\alpha_i,\alpha)$.  In particular, the dotted action of $W$ on 
$\check U^0$ specializes to the dotted (or translated) action of $W$ on $U(\gh)$ 
(see [H, Section 23.3]).   A similar argument shows that the dotted action of 
$W_{\Theta}$ on ${\cal C}[{\cal A}]$ specializes to the dotted action 
of $W_{\Theta}$ on $U(\ga)$ as defined in [He, Section 2].  Thus we have a 
commutative diagram
$$\matrix{{\cal C}[T_2]^{W.}&\rightarrow &U(\gh)^{W.}\cr
	     \ &\ &\ \cr
	     \downarrow&\ &\downarrow\cr
	     \ &\ &\ \cr
	      {\cal C}[{\cal A}]^{W_{\Theta}.}&\rightarrow 
	      &U(\ga)^{W_{\Theta}.}}\eqno{(5.8)}$$ 
where the maps are the same as in the commutative diagram (5.7).
Note that the top and bottom maps are still surjections.

\begin{theorem}
    Suppose $\gg,\gg^{\theta}$ is an irreducible symmetric pair of 
    type EIII, EIV, EVII, or EIX.   Then ${\cal P}_B(Z(\check 
    U))\neq {\cal C}[{\cal A}]^{W_{\Theta}.}$
    \end{theorem}
    
    \noindent{\bf Proof:}   Note that the image of the center 
    $Z(\check U)$ under the ordinary (quantum) Harish-Chandra map 
    ${\cal P}$ is 
    just ${\cal C}[T_2]^{W.}$ ([Jo, Lemma 7.1.17 and 7.1.25]).   Assume 
    that ${\cal P}_B(Z(\check 
    U))$ $= {\cal C}[{\cal A}]^{W_{\Theta}.}$.  Recall that the 
    restriction of ${\cal P}_B$ to $Z(\check U)$ agrees 
    ${\cal P}_B\circ {\cal P}$. 
    It follows that  the map 
    from ${\cal C}[T_2]^{W.}$ to ${\cal C}[{\cal A}]^{W_{\Theta}.}$ 
    in the commutative diagram (5.8) is surjective.   Hence the map 
    from $U(\gh)^{W.}$ to $U(\ga)^{W_{\Theta}.}$ is surjective.  
    Therefore the image of the center of $U(\gg)$ under the 
    classical Harish-Chandra map associated to $\gg,\gg^{\theta}$ is 
    equal to  the entire invariant ring $U(\ga)^{W_{\Theta}.}$.  This 
    contradicts [He].   $\Box$
    
    \medskip
    It should be noted that   surjectivity of   the map from $U(\gh)^{W.}$ to 
    $U(\ga)^{W_{\Theta}.}$ does not imply surjectivity of    the map from 
    ${\cal C}[T_2]^{W.}$ to ${\cal C}[{\cal A}]^{W_{\Theta}.}$.  The problem is that there are many subrings of 
    ${\cal C}[{\cal A}]^{W_{\Theta}.}$ which surject onto 
    $U(\ga)^{W_{\Theta}.}$.  (One such example is ${\cal 
    C}[t^2|t\in {\cal A}]^{W_{\Theta}.}$.)   The reader is referred 
    to [JL2, 6.13] for a similar situation.  In particular,    it is shown 
    in [JL2] that the 
    center $Z(U)$ of the ordinary quantized enveloping algebra $U$  
    specializes to the center $Z(\gg)$ of $U(\gg)$.  However, for 
    many simple Lie algebras $\gg$,  $Z(U)$ is not 
      a 
    polynomial ring while $Z(\gg)$ is always a polynomial ring. 
 
\section{The image of the center}

 In the classical case, the image of the center of $U(\gg)$ under the 
 Harish-Chandra map associated to the symmetric pair 
 $\gg,\gg^{\theta}$ is equal to the image of the whole invariant 
 subalgebra $U(\gg)^{\gg^{\theta}}$ in all but the  four exceptional cases
 EIII, EIV, 
EVII, and EIX.
 The purpose of  this section  is to establish 
 the quantum version of this result.   
In particular, we prove the following 
generalization of Theorem  5.4 and Theorem 5.5.

\begin{theorem} Let  $\gg,\gg^{\theta}$ be an 
irreducible symmetric pair.  Then for all $B\in {\cal B}$,  ${\cal P}_B(Z(\check 
U))={\cal 
C}[{\cal A}]^{W_{\Theta}.}$ if and only if $\gg,\gg^{\theta}$ is not 
of type EIII, EIV, EVII, and EIX.  Moreover, in these cases,
${\cal P}_B(\check U^B)={\cal P}_B(Z(\check 
U)).$
\end{theorem}

Let $B$ be a fixed coideal subalgebra in the set ${\cal 
B}$ associated to $\gg,\gg^{\theta}$.  Write $t$ for the rank of the restricted root system 
associated to $\gg,\gg^{\theta}$.
Note that when $\gg,\gg^{\theta}$ is of type  EIII, EIV, EVII, EIX or
$\widetilde{P^+(\pi)}=P^+(\Sigma)$, then the assertions of the theorem 
follow from Theorem 5.4 and Theorem 5.5. Thus we assume in this 
section that  
 $\gg$ contains a $\theta$ invariant Lie subalgebra $\gr$ 
such that $\gr,\gr^{\theta}$ is 
of type AII   and the rank of the restricted root system associated to 
$\gr, \gr^{\theta}$ is $t-1$.  We further assume that $t\geq 4$ 
and $\gg$ is of classical type. It should be noted that all
  irreducible symmetric pairs not covered in the theorems of Section 
5 satisfy these conditions. Thus, we prove Theorem 6.1 by showing that ${\cal P}_B(Z(\check 
U))={\cal 
C}[{\cal A}]^{W_{\Theta}.}$  under precisely these 
assumptions. 

The possible symmetric pairs that satisfy 
the   assumptions presented in the previous paragraph are given in (3.11)-(3.15). As in Section 3, we 
write $\eta_1,\dots, \eta_t$ for the fundamental restricted roots 
and $\mu_1,\dots, \mu_t$ for the simple restricted roots.  We also 
assume that
  $\eta_0=0$
and $\eta_{t+1}=0$.

Now $Z(\check U)$ is generated by $z_{2\omega_1},\dots, 
z_{2\omega_n}$ ([Jo, Section 7.3]).  Thus, Theorem 2.4 ensures that   
${\cal P}_B(Z(\check U))$  is generated by ${\cal P}_B(z_{2\omega_1}),\dots, 
{\cal P}_B(z_{2\omega_n})$.  On the other hand, ${\cal C}[{\cal 
A}]^{W_{\Theta}.}$ is generated by $\hat m(2\eta_i), 1\leq i\leq n$.   
Therefore, in order to show ${\cal P}_B(Z(\check U))$ is equal to 
${\cal C}[{\cal 
A}]^{W_{\Theta}.}$, it is sufficient to  express the $\hat 
m(2\eta_i), 
1\leq i\leq n$ as polynomials in   the ${\cal P}_B(z_{2\omega_i})$, 
$1\leq i\leq n$.   The first step  is to write the ${\cal P}_B(z_{2\omega_i})$
as a linear combination of the $\hat m(2\gamma)$ for $\gamma\in 
P^+(\Sigma)$.   The next  lemma which analyzes dominant integral weights of 
the form $\eta_r+\eta_s$ is critical in this process.

\begin{lemma}
 Assume that $\Sigma$ is of type $A_t$.
 Suppose that $r$ and $k$ are    positive integers so that $r\leq k$ 
 and $r+k\leq 
t+1$. The set of elements of  
$P^+(\Sigma)$ which are strictly less than $ \eta_r+\eta_k$ is 
$\{ \eta_{r-s}+\eta_{k+s}|1\leq s\leq r\}$. 
\end{lemma}

\noindent
{\bf Proof:} 
It is straightforward to check that
$$\eta_r+\eta_k-\sum_{i=1}^s\sum_{j=r-i+1}^{k+i-1}\mu_j=
\eta_{r-s}+\eta_{k+s}\eqno(6.1)$$
for $1\leq s\leq r$.
Hence $\{\eta_{r-s}+\eta_{k+s}|1\leq j\leq r\}$ is a 
subset of the intersection of $\{\eta_{r}+\eta_k- \gamma|\gamma\in 
Q^+(\Sigma)\setminus\{0\}\}$ with 
$P^+(\Sigma)$.

Suppose that $i$ satisfies $1\leq i\leq t$. By [H, Section 13.2, Table 
1], the coefficient of $\mu_t$ in $\eta_i$ written 
as a linear combination of simple restricted roots $\mu_1,\dots, 
\mu_t$ is 
$i/(t+1)$.  Now suppose that $\gamma $ is a nonzero element of $ Q^+(\Sigma)$ such that $\eta_r+\eta_k-\gamma\in 
P^+(\Sigma)$.  
It follows that  the coefficient of $\mu_t$ in $ \eta_r+\eta_k$ written as a linear 
combination of the simple restricted roots   is 
$(r+k)/{(t+1)}$.  On the other hand, the fact that 
$\gamma\in Q^+(\Sigma)$ implies that the 
coefficient of $\mu_t$ in $\gamma$ is a nonnegative integer.   Since 
$  \eta_r+\eta_k-\gamma$ is in $P^+(\Sigma)$, the coefficient of $\mu_t$ in 
$ \eta_r+\eta_k-\gamma$ must be nonnegative.  Note that 
$0<(r+k)/(t+1)\leq 1$. It 
follows that the coefficient of $\mu_t$ in $\eta_r+\eta_k-\gamma$ 
written as a linear combination of the simple restricted roots 
must either be
$(r+k)/(t+1)$ or $0$.  This latter case can only occur if 
$r+k=t+1$.   Moreover, the only linear combination of the 
$\mu_1,\mu_2,\dots, \mu_{t-1}$ contained in   $P^+(\Sigma)$ is $0$.   Hence 
this latter case happens only when 
$\eta_r+\eta_k-\gamma=0=\eta_0+\eta_{t+1}$.    

Now assume that the coefficient of $\mu_t$ in $\eta_t+\eta_k-\gamma$ 
is $(r+k)/(t+1)$.
Write $\eta_r+\eta_k-\gamma=\sum_{i=1}^tm_i\eta_i$ where the $m_i$ are 
nonnegative integers. Using the previous paragraph, a comparison of the coefficient of $\mu_t$ in 
both sides of this equation yields 
$$\sum_i{{m_ii}\over{ (t+1)}}={{(r+k)}\over{(t+1)}}.\eqno(6.2)$$
We do a similar comparison using the restricted root $\mu_1$. 
By [H, Section 13.2, Table 1], the coefficient of $\mu_1$ in $ \eta_i$ is $(t-i+1)/(t+1)$.
Hence the coefficient of $\mu_1$ in $\sum_im_i\eta_i$ written as a 
linear combination of the $\mu_i$ is 
$$ 
\sum_i{{m_i(t-i+1)}\over{(t+1)}}=\sum_im_i-\sum_i{{m_ii}\over{(t+1)}}=\sum_im_i-{{(r+k)}\over{(t+1)}}$$
On the other hand, the coefficient of $\mu_1$ in $\eta_r+\eta_k$ is 
$${{(t-r+1)}\over{(t+1)}}+{{(t-k+1)}\over{(t+1)}}=2-{{(r+k)}\over{(t+1)}}.$$  
The fact that  $\eta_r+\eta_k> \sum_im_i\eta_i$ ensures that
$$0\leq \sum_im_i-{{(r+k)}\over{(t+1)}}< 2-{{(r+k)}\over{(t+1)}}.$$  

Since the  $m_i$ are 
nonnegative integers, there are two possibilities.   The first is that 
there exists an $i$ with $1\leq i\leq t$, 
  $m_i=1$, and   $m_j=0$ for $j\neq i$.  By (6.2),
 the coefficient $\mu_t$ in 
$ \eta_{r}+\eta_{k}-\gamma$ is $(r+k)/(t+1)$.  It follows from [H, 
Section 13.2, Table 1] that 
$i=r+k$. Thus $ \eta_{r}+\eta_{k}-\gamma=\eta_{ r+k}=\eta_0+\eta_{r+k}$.

In the second case, there exists
$i$ and $j$ with $1\leq i,j\leq t$,  $m_i=m_j=1$, and  $m_k=0$ for 
$k\notin\{i,j\}$.   Since the coefficient of $\mu_t$ in $
\eta_{r}+\eta_{k}-\gamma$ is $(r+k)/(t+1)$, we must have that 
$\eta_{r}+\eta_{k}-\gamma=\eta_{i}+\eta_{j}$  with $i+j= k+r$. If $i\geq r$, 
then (6.1) implies that $\eta_i+\eta_{j}\geq\eta_r+\eta_{k}$.   
Thus $i< r$.  In particular,  we can write $i=r-s$ and $j=k+s$ for some 
$1\leq
s\leq  r$.      $\Box$

\medskip
Consider the case when $\gg,\gg^{\theta}$ is an irreducible 
symmetric pair of type AII.   It follows that the restricted root 
system is of type $A_t$ and the  rank $n$ of $\gg$, 
  is 
equal to    
$2t+1$.  For this particular irreducible symmetric pair,
it is straightforward to check that the only elements of $Q(\pi)$ fixed by $\Theta$ are elements of 
$Q(\pi_{\Theta})$. It 
follows that for each $\alpha\in Q(\pi)$, 
$\tilde\alpha=0$ if and only if $\alpha\in Q(\pi_{\Theta})$. 
  
Given a  root $\alpha$ in the root system $\Delta$ generated by $\pi$, let $s_{\alpha}$ denote the 
reflection in $W$ associated to $\alpha$.
The next lemma determines the number of weights in 
the $W$ orbit of $\omega_{i+k}$ whose restriction is equal to 
$\eta_{i-s}+\eta_{k+s}$ where $k=i$ or $k=i+1$.

\begin{lemma} Assume that $\gg,\gg^{\theta}$ is an irreducible 
symmetric pair of type AII. Suppose $i$ is an integer such that 
$1\leq i\leq (t+1)/2$. For all integers  $s$ such that $0\leq 
s\leq i$, the number of elements in 
the set $$\{\beta\in W\omega_{2i}|
\tilde\beta=\eta_{i-s}+\eta_{i+s}\}$$ equals $2^s$.
Similarly, suppose $0\leq i\leq t/2$. Then for all $s$ such that 
$0\leq s\leq i$,  the number of of elements in 
the set $$\{\beta\in 
W\omega_{ 2i+1}|\tilde\beta=\eta_{i-s}+\eta_{i+1+s}\}$$ equals $2^{s+1}$.
\end{lemma}

\noindent
{\bf Proof:}   Fix $i$ such that $1\leq i\leq (t+1)/2$. Set 
$\gamma_0=\omega_{2i}$.   For each $s$ 
such that $1\leq s\leq i$, set 
$$\gamma_s=\omega_{2i-2s+1}-\omega_{2i-2s+2}+\omega_{2i-2s+3}+\cdots 
-\omega_{2i-2s+(4s-2)}+\omega_{2i-2s +(4s-1)}.$$
Recall that 
$\mu_j=\tilde\alpha_{2j}$ for $1\leq j\leq 
t$.  Hence by (3.10), if $2j+1$ is not an element of the 
set 
$\{2i-2s+1, 2i-2s+2, 2i+2s-2, 2i+2s-1\}$, then $(\gamma_s,\mu_j)=0$.  Since $j$ is an integer, 
$2j+1$ cannot equal $2i-2s+2$ or $2i+2s-2$.  
Temporarily set $\mu_0=\mu_{t+1}=0$. If $2j+1=2i-2s+1$, then $j=i-s$ and  
$(\gamma_s,\mu_j)=(\omega_{2i-s},\mu_j)=(\eta_j,\mu_j)$.   Similarly,
if $2j-1=2i+2s-1$ then $j=i+s$ and  
$(\gamma_s,\mu_j)=(\omega_{2i+s},\mu_j)=(\eta_j,\mu_j)$.   Thus 
$\tilde\gamma_s=\eta_{i-s}+\eta_{i+s}$. (Note that 
these computations work  even in the special case  when $i=s$.)   

We claim that $\gamma_s\in W\omega_{2i}$ for $0\leq 
s\leq i$.  Now $\gamma_0=\omega_{2i}$ and so $\gamma_0\in 
W\omega_{2i}.$   Also, $(\gamma_0,\alpha_{2i})=1$ and in particular, 
$s_{\alpha_{2i}}\gamma_0=\gamma_0-\alpha_{2i}=\omega_{2i-1}-\omega_{2i}+\omega_{2i+1}=\gamma_1$.
Hence $\gamma_1\in W\omega_{2i}$.   Now assume that $\gamma_s\in 
W\omega_{2i}$.   It is straightforward to check that 
$$\gamma_{s+1}=\gamma_s-\alpha$$ where
$$\alpha= \alpha_{2i-2s-1}+\alpha_{2i-2s}+\alpha_{2i-2s+1}+\cdots 
+\alpha_{2i+2s+1}.$$  Note that $\alpha$ is a positive root in 
$\Delta$. Moreover, $(\gamma_s,\alpha)=1$.  Hence
$s_{\alpha}\gamma_s=\gamma_{s+1}$.  Therefore, the claim follows by 
induction on $s$.

Note that the set $$\{\beta\in W\omega_{2i}|
\tilde\beta=\eta_{i-s}+\eta_{i+s}\}\eqno{(6.3)}$$  equals the set 
$$\{\beta\in W\omega_{2i}|
\tilde\beta=\tilde\gamma_s\}.$$   By the discussion preceding the 
theorem,
$\tilde\beta=\tilde\gamma_s$ if and only if $\beta-\gamma_s\in 
Q(\pi_{\Theta})$. Recall that $W'$ is the Weyl group associated to 
the root system generated by $\pi_{\Theta}$.   Since $\beta\in W\gamma_s$ and 
$\gamma_s-\beta\in Q(\pi_{\Theta})$, it follows that  
$\beta\in W'\gamma_s$.
Hence the number of elements in (6.3) is equal to 
$|W'|/|{\rm Stab}_{W'}\gamma_s|$. 

Set $s_i=s_{\alpha_i}$ for $\alpha_i\in \pi$.  Note that $W'$ is generated by 
$\{s_{{2r+1}}|0\leq r\leq t\}$.   Morever, each $s_{{2j+1}}$ has order 
$2$ while $s_{2j+1}s_{2r+1}=s_{2r+1}s_{2j+1}$ for $r\neq j$.   Hence 
$W'$ is isomorphic to $t+1$ copies of ${\bf Z}_2$.   Thus 
$|W'|=2^{t+1}$.   On the other hand, ${\rm Stab}_{W'}\gamma_s=\{s_{2r+1}| 
r\neq 2i+2j+1$ for $-s\leq j\leq s-1\}.$ It follows that 
$ {\rm Stab}_{W'}\gamma_s $ is isomorphic to $t+1-(2s)$ copies 
of ${\bf Z}_2$.  Hence $|{\rm Stab}_{W'}\gamma_s|=2^{t+1-2s}$ and 
$|W'|/|{\rm Stab}_{W'}\gamma_s|=2^{2s}.$  This completes the proof of the 
first assertion. The second assertion follows in a similar 
fashion.$\Box$

\medskip
We continue the analysis of the case when $\gg,\gg^{\theta}$ is of type  AII.
Note that  $\gg,\gg^{\theta}$   of type  AII  implies
that $\gg$ is of type  A$_n$.   Hence $\omega_{i}$ is a minuscule 
fundamental 
weight for $1\leq i\leq n$. 
(See [M] or [L5, Section 7] for a definition of minuscule weight.)  In particular, there does not exist $\gamma\in P^+(\pi)$ 
such that $\omega_{i}>\gamma$.  By Lemma 5.1, there exists a nonzero 
scalars $a_i$ such 
that ${\cal P}(a_iz_{2\omega_{i}})$ is just   the sum 
 $$\sum_{\mu\in W\omega_{i}}q^{(\rho,2\mu)}\tau(2\mu).$$   
 By the discussion following the proof of Lemma 5.1, ${\cal P}_B(a_iz_{2\omega_{i}})$ is a linear 
 combination of terms of the form $\tau(2\tilde\mu)$ with $\tilde\mu\in 
 P(\Sigma)$.   Moreover,  the coefficients are Laurent polynomials in $q$ with 
 nonnegative integer coefficients and the coefficient of 
 $\tau(2\tilde\omega_i)$ is nonzero.   Lemma 5.3 ensures that we can 
 also write ${\cal P}_B(a_iz_{2\omega_{i}})$ as a linear 
 combination of terms of the form $\hat m(2\tilde\mu)$ with $\tilde\mu\in 
 P^+(\Sigma)$. The next lemma gives 
 the evaluation  at $q=1$ of 
 the  coefficients  of the $\hat m(2\tilde\mu)$  in these sums for some 
 of the ${\cal P}_B(a_iz_{2\omega_{i}})$.

 \begin{lemma}   Assume that $\gg,\gg^{\theta}$ is an irreducible 
 symmetric pair of type AII.      Then  
 for each integer $i$ such that $2\leq 2i\leq t+1$, there exist Laurent 
 polynomials $f_{2s}(q), 0\leq s\leq i$, such that
 ${\cal P}_B(z_{2\omega_{2i}})$ is a nonzero scalar multiple 
  of 
  $$\sum_{0\leq s\leq i}f_{2s}(q) 
 \hat m(2\eta_{i-s}+2\eta_{i+s}).$$
Similarly, for each integer $i$ such that $1\leq 2i+1\leq t+1$, there exist Laurent 
 polynomials $f_{2s+1}(q), 0\leq s\leq i$, such that  ${\cal P}_B(z_{2\omega_{2i+1}})$ is a nonzero scalar 
 multiple of 
 $$\sum_{0\leq s\leq i}f_{2s+1}(q) 
 \hat m(2\eta_{i-s}+2\eta_{i+1+s}). $$  Moreover, the coefficients of 
 powers of $q$ in each $f_j(q)$ are 
 nonnegative integers and $f_j(1)=2^j$ for all $0\leq j\leq 2i+1$.
\end{lemma}

\noindent
{\bf Proof:} 
 Let $i$ be a positive integer such that $i\leq (t+1)/2$.  By  Lemma 6.2,
the only restricted dominant integral weights less than or equal to 
 $\tilde\omega_{2i}$ are of the form $\eta_{i-s}+\eta_{i+s}$ for 
 $0\leq s\leq i$.   Hence, by the discussion preceding the lemma, for each $s$ such that $0\leq s\leq i$, 
 there exists a Laurent polynomial $f_{2s}(q)$ in $q$ such that
 $${\cal P}_B(a_iz_{2\omega_{2i}})=\sum_{0\leq s\leq i}f_{2s}(q) 
 \hat m(2\eta_{i-s}+2\eta_{i+s}).$$  
 Recall that each $W_{\Theta}$ orbit of a single element in 
 $P(\Sigma)$ contains a single dominant integral weight in 
 $P^+(\Sigma)$.   Hence if we expand out ${\cal P}_B(a_iz_{2\omega_{2i}})$ as a sum of terms of the form 
 $\tau(2\tilde\mu)$ with $\tilde\mu\in P(\Sigma)$, it follows that the 
 coefficient of $\tau(2\eta_{i-s}+2\eta_{i+s})$ is $f_{2s}(q)$.  On 
 the other hand, 
 this coefficient can be computed using the expression of  
    ${\cal P}(az_{2\omega_{2i}})$ in the   paragraph preceding the 
    lemma.  In particular, 
    $f_{2s}(q)$ is just the sum  
    $$\sum_{\{\tilde\gamma\in W\omega_{2i}|\tilde\gamma 
    = \eta_{i-s}+ \eta_{i+s}\}}q^{(\rho,2\tilde\gamma)}.$$ This shows 
    that each $f_{2s}(q)$ is a Laurent polynomial in $q$ with 
    nonnegative coefficients.   Moreover,  
    $f_{2s}(1)$ is just the number of elements $\gamma$ in $W\omega_{2i}$ such 
    that $\tilde\gamma= \eta_{i-s}+ \eta_{i+s}$.  By Lemma 6.3, we have 
    $f_{2s}(1)=2^{2s}$.  This completes the proof 
   for the central element $z_{2\omega_{2i}}$. The proof for  
   $z_{2\omega_{2i+1}}$ follows in a similar 
    fashion. $\Box$
     
    \medskip

Note that the map which sends
$q^{(\rho,2\eta_i)}\tau(2\eta_i) $ to 
$z^{\eta_i}$ for each $1\leq i\leq t$ defines an algebra isomorphism 
from 
${\cal C}[{\cal A}]$ onto ${\cal 
C}[P(\Sigma)]$.    For each $\lambda\in 
P(\Sigma)$, set $$m(\lambda)=\sum_{\gamma\in 
W_{\Theta}\lambda}z^{\gamma}.$$   The image of  $\hat m(2\lambda)$ 
under this isomorphism is just 
  $m(\lambda)$.  In the next few lemmas, the 
computations are done using $m(\lambda)$ instead of $\hat m(2\lambda)$ 
in order to make the notation easier.

\begin{lemma} Assume that $\gg,\gg^{\theta}$ is an irreducible 
symmetric pair of type AII. 
    Let $r$ and $k$ be two positive integers such that $r\leq k$ and 
    $r+k= t+1$.
      Then for all $0\leq j\leq r$, the coefficient of $1$ in    
    $\hat m(2{\omega_{r-j}) \hat m(2\omega_{k+j}})$ written as a linear combination 
    of the set $\{\hat m(2\lambda)|\lambda\in P^+(\Sigma)\}$ is 
    ${{t+1}\choose{r-j}}$.
       \end{lemma}
    
   \noindent{\bf Proof:}  Using the isomorphism described above, we replace each
   $\hat m(2\lambda)$ with  $  m(\lambda)$ in the proof of the lemma.
   Fix $j$ such that $0\leq j\leq r$.  Note that if $j=r$,
   then $ m(\eta_{r-j}) m(\eta_{k+j})= m(\eta_{0}) m(\eta_{t+1})=1$.   
   Since ${{t+1}\choose {0}}=1$, the lemma follows in this case.  
   Hence we may assume that $j<r$.
   
   Note that 
   $$ m(\eta_{r-j})\in
   z^{\eta_{r-j}}+\sum_{\gamma<\eta_{r-j}}{\bf 
   N}z^{\gamma}
   $$
   and 
   $$ m(\eta_{k+j})\in z^{\eta_{k+j}}+\sum_{\gamma<\eta_{k+j}}{\bf 
   N}
  z^{\gamma}.$$
   Hence $$m(\eta_{r-j})m(\eta_{k+j})\in z^{\eta_{r-j}+\eta_{k+j}}
   +\sum_{\gamma<\eta_{r-j}+\eta_{k+j}}{\bf 
   N}z^{\gamma}.$$ 
    Now $m(\eta_{r-j})m(\eta_{k+j})$ is $W_{\Theta}$ invariant.   
    Hence Lemma 6.2 
   implies that there exists a nonnegative integer $a$ such that 
   $$\eqalign{m(\eta_{r-j})m(\eta_{k+j})\in\ & 
   m(\eta_{r-j}+\eta_{k+j})+am(0)\cr&+\sum_{1\leq b< r-j} {\bf 
   N}m(\eta_{r-j-b}+\eta_{k+j+b}).\cr}$$  Furthermore, $m(0)$ is just 
   $z^0=1$ and  $a$ is  the coefficient of $1$ in 
   $m(\eta_{r-j})m(\eta_{k+j})$. 
   
  Since the root system $\Sigma$ is of type $A_{t}$ and $r+k=t+1$, 
  $\Sigma$  
   admits a diagram automorphism $d$ which sends $\mu_{r-j}$ to 
   $\mu_{k+j}$ for $0\leq j\leq r-1$.   Let $w'_0$ denote the 
   longest element of the Weyl group $W_{\Theta}$ and note that $w'_0=-d$.  We 
   have the following equality of sets $$\{w\mu_{r-j}|w\in 
   W_{\Theta}\}=\{ww'_0\mu_{r-j}|w\in W_{\Theta}\}=\{-w\mu_{k+j}|w\in 
   W_{\Theta}\}.$$
   Thus if we write $$m(\mu_{r-j})=\sum_{1\leq i\leq 
   s}z^{\gamma_i},$$ then $$m(\mu_{k+j})=\sum_{1\leq i\leq 
   s}z^{-\gamma_i}.$$ It follows that $s=a$.
   Moreover, $s$ is just the number of elements 
   in the orbit $W_{\Theta} \mu_{r-j}$.   Hence $s=|W_{\Theta}|/|{\rm 
   Stab}_{W_{\Theta}} \mu_{r-j}|$.  Now ${\rm Stab}_{W_{\Theta}}\mu_{r-j}$ 
   is just the 
   subgroup of $W_{\Theta}$ generated by the reflections corresponding 
   to the 
   simple roots $\mu_k$ with $k\neq {r-j}$.   In particular, 
   ${\rm Stab}_{W_{\Theta}}\mu_{r-j}$ is isomorphic to the direct product of 
   two groups $W_1\times W_2$ where  $W_1$  is the Weyl group 
   associated to a root system of type $A_{r-j-1}$ and $W_2$ is the 
   Wey group associated to a root system of type $A_{t-r+j}$.
   Hence, by [H, Section 12.2, Table 1] it follows that $$a=|W_{\Theta}| |{\rm Stab}_{W_{\Theta}}\mu_{r-j}|^{-1}
   ={{(t+1)!}\over{(r-j)!(t+1-r+j)!}}
   = {{t+1}\choose{ r-j}}.\Box$$

   \medskip  
 Suppose that  $\Sigma'_s$ is a subset   of  the set of  simple restricted roots $\{\mu_1,\dots, 
  \mu_t\}$.  We define here notation associated to $\Sigma'_s$ which 
  is used in the rest of the section for different choices of 
  subsets $\Sigma'_s$ of the set of simple restricted roots. Let $\Sigma'$ 
  be the  sub root system of $\Sigma$  generated by $\Sigma_s'$. 
  Let $W'_{\Theta}$ denote the Weyl group 
  of  $\Sigma'$.  For each 
  $\lambda\in P(\Sigma)$, let $\lambda'$ be the element of $P^+(\Sigma')$ 
  such that   
  $(\lambda-\lambda',\gamma)=0$ for all $\gamma\in Q(\Sigma')$. 
  Set $${\cal N}=\sum_{\mu_i\in \Sigma_s\setminus 
   \Sigma'_s}\sum_{\gamma\in Q^+(\Sigma)}{\bf N} 
   z^{-\mu_i-\gamma}.$$
   Given $\lambda\in P^+(\Sigma)$, 
   set $$m'(\lambda)= \sum_{\gamma\in 
  W_{\Theta}'\lambda} z^{\gamma}.$$  Note that   $$m(\lambda)\in 
  z^{\lambda-\lambda'}
  m'(\lambda')+z^{\lambda}{\cal N}.\eqno{(6.4)}$$  Moreover, 
  $$m(\lambda)m(\beta)\in 
z^{\lambda-\lambda'+\beta-\beta'}m'(\lambda')m'(\beta')
 +z^{\lambda+\beta}{\cal N}.\eqno{(6.5)} $$

 Consider for the moment the special case when $\Sigma_s$ is the 
 empty set.   Expression (6.5) becomes 
$$m(\lambda)m(\beta)\in 
z^{\lambda+\beta} 
 +z^{\lambda+\beta}\ \sum_{\gamma\in Q^+(\Sigma)\setminus\{0\}}{\bf N} 
   z^{-\gamma}. $$  Now suppose that $\lambda=\eta_r$ and 
   $\beta=\eta_k$.   It follows from Lemma 6.2 that 
   $m(\eta_r)m(\eta_k)$ can be written as a linear combination of 
   elements in the set $\{m(\eta_{r-s}+\eta_{k+s})|\ 0\leq s\leq r\}$.

   \begin{lemma}  
   Assume that $\gg,\gg^{\theta}$ is an irreducible symmetric pair of 
   type AII.
    Let $r$ and $k$ be two positive integers such that $r\leq k$ and 
    $r+k\leq  t+1$.
    Then  
    $$\hat m(2{\eta_{r})\hat m(2\eta_{k}})=\sum_{0\leq s\leq r}{{k-r+2s}\choose{s}}\hat m(2\eta_{r-s}+2\eta_{k+s}).$$ 
       \end{lemma}
       
   \noindent
   {\bf Proof:} Once again, we use $m(\lambda)$ instead of $\hat 
   m(2\lambda)$ in proving this lemma. 
   First consider the special case when $r=0$ and $k=t+1$.   Note 
   that this forces $s=0$.   Moreover, $m({\eta_{0})m(\eta_{n+1}})
   =m(0)m(0)=1$ while  ${{t+1}\choose{0}}
   =1$.   So the lemma follows in this case.   Now 
   assume that $r$ and $s$ are   chosen so that $0\leq s\leq r$ and 
   either $r>0$ or $k<t+1$.

   Set $\Sigma'_s=\{\mu_{i}|\ r-s+1\leq i\leq k+s-1\}$ and note 
   that  
  $\Sigma_s'$ is the set of simple roots for a root system of type 
 $A_{k-r+2s-1}$.  
 For each $i$, the weight $\eta_i'$  is just the fundamental weight corresponding to 
 the simple root 
 $\mu_i$ in the weight lattice associated to the root system 
 $\Sigma'_s$.  Checking the
 formulas for the fundamental weights in [H, Section 3.2, Table 1] yields that 
 $\eta'_{r}+\eta'_{k}\in Q^+(\Sigma')$.  Since $Q^+(\Sigma')$ is  a 
 subset of $Q^+(\Sigma)$, it further follows that $\eta'_{r}+\eta'_{k}\in Q^+(\Sigma)$.
 Now $(\mu_i, \mu_b)\leq 0$ for all $\mu_i\in \Sigma'$ and 
 $\mu_b\in \Sigma\setminus \Sigma'$.   Thus 
 $(-\eta'_{r}-\eta'_{k}, \mu_b)$ is a nonnegative integer
 for all $\mu_b\in 
 \Sigma\setminus \Sigma'$. Moreover, since $\Sigma$ is 
 a root system of type $A_t$,  $\mu_b\in 
 \Sigma\setminus \Sigma'$ implies that $(-\eta'_{r}-\eta'_{k}, \mu_b)=0$ 
 unless $b=r-s$ or $b=k+s$.  
 
 Now consider the weight 
 $\gamma=\eta_{r}+\eta_{k}-\eta'_{r}-\eta'_{k}$.  
 The above discussion  ensures that $(\gamma, \mu_b)=0$ for all 
 $\mu_b\in \Sigma\setminus\{\mu_{r-s},\mu_{k+s}\}$
 while $(\gamma,\mu_b)\geq 0$ for  
 $\mu_b\in \{\mu_{r-s},\mu_{k+s}\}$.    Thus $\gamma$ is linear 
 combination of the fundamental weights $\eta_{r-s}$ and 
 $\eta_{k+s}$ with nonnegative integer coefficients.  
   By the previous 
 paragraph,
  $(\eta_{r}+\eta_{k})-\gamma=\eta'_{r}+\eta'_{k}$ is in $Q^+(\Sigma)$.  
  Hence it follows from  
 Lemma 6.2 that $\gamma$ must equal $\eta_{r-s}+
  \eta_{k+s}$.  In 
 particular, we have 
 $$\eta_{r}+\eta_{k}-\eta'_{r}-\eta'_{k}=\eta_{r-s}+
  \eta_{k+s}.\eqno{(6.6)}$$ 
 Thus, applying (6.5) to $m(\lambda)m(\beta)$  yields 
  $$m(\eta_{r})m(\eta_{k})\in 
   z^{\eta_{r-s}+
  \eta_{k+s}}m'(\eta'_{r})m'(\eta'_{k})+
  z^{\eta_{r}+\eta_{k}}{\cal N}.$$

Consider the following two ways to expand the product 
$m(\eta_{r})m(\eta_{k})$.  The first  is  as a 
linear combination of terms of the form $z^{\beta}$ for $\beta\in 
P(\Sigma)$ while the second  is as a 
linear combination of terms of the form $m(\lambda)$ for $\lambda\in 
P^+(\Sigma)$.     Note that the coefficient of $z^{\eta_{r-s}+
  \eta_{k+s}}$ in  $m(\eta_{r})m(\eta_{k})$ written the first 
  way is the same as the coefficient of 
  $m(\eta_{r-s}+\eta_{k+s})$ in $m(\eta_{r})m(\eta_{k})$  written the 
  second way.  
  By (6.6),  $z^{\eta_{r-s}+
  \eta_{k+s}}$ is not an element of  
  $z^{\eta_{r}+\eta_{k}}{\cal N}$.   Hence  the coefficient of $z^{\eta_{r-s}+
  \eta_{k+s}}$ in $m(\eta_{r})m(\eta_{k})$ equals the coefficient of 
 $z^{\eta_{r-s}+
  \eta_{k+s}}$ in  $z^{\eta_{r-s}+
  \eta_{k+s}}m'(\eta'_{r})m'(\eta'_{k})$.   But this is 
  the same as the coefficient of $1$ in $m'(\eta'_{r})m'(\eta'_{k})$.
 Recall that the first simple root in $\Sigma_s'$ is $\mu_{r-s+1}$ and 
 so $\mu_{r}$ is the $s^{th}$ simple root.  Also, there are 
 $k-r+2s-1$ simple roots in $\Sigma_s'$.  
Thus, by Lemma 6.5, the coefficient of $z^{\eta_{r-s}+
  \eta_{k+s}}$ in $z^{\eta_{r-s}+
  \eta_{k+s}}m'(\eta'_{r})m'(\eta'_{k})$ is 
  $${{k-r+2s}\choose{s 
  }}   \Box$$
  
 \medskip
 Assume that $\gg,\gg^{\theta}$ is of type AII.  Given an integer $l$ such that $1\leq 2l\leq 
 t+1$, set 
   $$\hat M_{2l}  =\sum_{0\leq i\leq 
   r}2^{2i} \hat m(2\eta_{l-i}+2\eta_{l+i}).\eqno{(6.7)}$$
  Similarly, given an integer $l$ satisfying $1\leq 2l+1\leq t+1$, set $$\hat M_{2l+1}  =\sum_{0\leq i\leq 
   r}2^{1+2i}\hat  m(2\eta_{l-i}+2\eta_{l+1+i}).\eqno{(6.8)}$$  By Lemma 
   6.4,   $\hat 
   M_{k}$ can be obtained from    ${\cal P}_B(z_{2\omega_k})$ by 
   writing the latter as a linear combination of  
  the $\hat m(2\lambda)$, for $\lambda\in 
  P^+(\Sigma)$, and evaluating each coefficient at  $q=1$. Let $M_{k}$ be the image of $\hat 
   M_k$ in ${\cal C}[P(\Sigma)]$ obtained using the isomorphism   described 
   before Lemma 6.5.

   \begin{lemma} Assume that $\Sigma$ is of type $A_t$.
       \begin{enumerate}
	   \item[(i)] Let $l$ be a
         positive integer such that $1\leq 2l\leq t+1$.
   Then  $\hat m(\eta_{2l})$ is in the   span of
   the set $\{{\cal P}_B(z_{2\omega_{2l}})\}\cup\{  
    \hat m(2{\eta_{l-j})\hat m(2\eta_{l+j}})|\ 0\leq j<l\}.$
     \item[(ii)] Let $l$ be a
         positive integer such that $1\leq 2l+1\leq t+1$.
   Then  $\hat m(2\eta_{2l+1})$ is in the  span of
   the set $\{{\cal P}_B(z_{2\omega_{2l+1}})\}\cup\{  
    \hat m(2{\eta_{l-j})\hat m(2\eta_{l+1+j}})|\ 0\leq j<l\}.$
    \end{enumerate}
        \end{lemma}
	
\noindent
{\bf Proof:} We prove (i).    Assertion (ii) follows using a similar 
argument.   We first prove that $  m(\eta_{2l})$ is in the ${\bf Q}$   span of
   the set $\{M_{2l}\}\cup\{  
    m({\eta_{l-j}) m(\eta_{l+j}})|\ 0\leq j<l\}.$
    Recall  that 
  $${{2i}\choose{r 
  }}={{2i}\choose{2i-r 
  }}$$ for all $0\leq r\leq 2i$.   Hence
  $$2^{2i}=(1+1)^{2i}=\sum_{0\leq r\leq 2i}{{2i}\choose{r 
  }}={{2i}\choose{i 
  }}+2\sum_{0\leq r\leq i-1}{{2i}\choose{r 
  }}.$$   By Lemma 6.6, it follows that  
  $$\eqalign{
  &m(\eta_l)m(\eta_l)+\sum_{1\leq j\leq l-1} 
  2m({\eta_{l-j})m(\eta_{l+j}})\cr&=
  \sum_{0\leq s\leq l} {{2s}\choose{s 
  }}m(\eta_{l-s}+\eta_{l+s})+\sum_{1\leq j\leq l-1} 
  \sum_{0\leq s\leq l-j} 2{{2j+2s}\choose{s 
  }}m(\eta_{l-j-s}+\eta_{l+j+s})\cr
 &=\sum_{0\leq s\leq l} {{2s}\choose{s 
  }}m(\eta_{l-s}+\eta_{l+s})+\sum_{1\leq j\leq l-1} 
  \sum_{j\leq i\leq l} 2{{2i}\choose{i-j 
  }}m(\eta_{l-i}+\eta_{l+i}).
  \cr}$$   One checks that   the subset  $\{(j,i)|
  1\leq j\leq l-1$ and $j\leq i\leq l\}\cup \{(l,l)\}$ of ${\bf Z}\times {\bf Z}$ 
  is equal to the subset $\{(j,i)| 1\leq i\leq l$ and $1\leq j\leq 
  i\}.$ Hence   
  $$\eqalign{
  &m(\eta_l)m(\eta_l)+\sum_{1\leq j\leq l-1} 
  2m({\eta_{l-j})m(\eta_{l+j}})\cr &=\sum_{0\leq s\leq l} {{2s}\choose{s 
  }}m(\eta_{l-s}+\eta_{l+s})+\sum_{1\leq i\leq l} 
  \sum_{1\leq j\leq i} 2{{2i}\choose{i-j 
  }}m(\eta_{l-i}+\eta_{l+i})-2{{2l}\choose{0 
  }}m(\eta_{2l})\cr
   &=\sum_{0\leq s\leq l} 
 \left({{2s}\choose{s 
  }}+\sum_{1\leq j\leq s} 2{{2s}\choose{s-j 
  }}\right)m(\eta_{l-s}+\eta_{l+s})-2m(\eta_{2l})\cr
  &=\sum_{0\leq s\leq l} 
 2^{2s}m(\eta_{l-s}+\eta_{l+s})-2m(\eta_{2l})
 =M_{2l}-2m(\eta_{2l})\cr}$$

Define $X$   by  
  $$X={\cal P}_B(z_{2\omega_{2l}})-\hat 
  m(2\eta_l)\hat m(2\eta_l)-\sum_{1\leq j\leq l-1} 
  2\hat m({2\eta_{l-j})\hat m(2\eta_{l+j}}).$$ 
 By Lemma 6.4  and Lemma 6.6 we can write $X-2\hat m(2\eta_{2l})$ as a linear combination
 of elements in the set $\{\hat 
  m(2\eta_{l-i}+2\eta_{l+i})|0\leq i\leq l\}$. Set $N$ equal to the ${\bf Q}$ vector space $\sum_{0\leq i\leq l}{\bf Q}\hat 
  m(\eta_{l-i}+\eta_{l+i}).$ 
The above 
  computations ensure that  
 the coefficient of each $\hat 
  m(2\eta_{l-i}+2\eta_{l+i})$ in $X$ is an element of $(q-1){\bf 
  Q}[q,q^{-1}]N$.  Now Lemma 6.6
  implies that the set $$\{\hat m(2\eta_{l-i})\hat m(2\eta_{l+j})+{\bf 
  Q}\hat m(2\eta_{2l})|0\leq i<l\} $$ is a basis for the ${\bf 
  Q}$ vector space $N/({\bf Q}\hat m(2\eta_{2l}))$. Hence there exist
  Laurent polynomials $g_0,\dots, g_{l}$ in ${\bf Q}[q,q^{-1}]$ so 
  that 
   $$X-2\hat m(2\eta_{2l})=\sum_{0\leq i\leq l-1}(q-1)g_i\hat 
   m(2\eta_{l-i})\hat m(2\eta_{l+j})+(q-1)
  g_l\hat m(2\eta_{2l}).$$
  Thus $$X-\sum_{0\leq i\leq l-1}(q-1)g_i\hat 
   m(2\eta_{l-i})\hat m(2\eta_{l+j})=(2+(q-1)
  g_l)\hat m(2\eta_{2l}).$$
  The lemma now follows from the fact that $2+(q-1)g_l$ cannot be zero.
    $\Box$
   
   \medskip
   
   We use  Lemma 6.7 to show that a certain set generates 
   ${\cal C}[{\cal A}]^{W_{\Theta}.}$ when $\gg,\gg^{\theta}$ is of type 
   AII.
   
   \begin{theorem}  Assume that $\gg,\gg^{\theta}$ is of type  AII.  Then
        ${\cal P}_Z(\check U))=
       {\cal C}[{\cal A}]^{W_{\Theta}.}$.
       \end{theorem}
       
       \noindent
       {\bf Proof:} By Theorem 2.6 and Lemma 5.3, we have that ${\cal 
       P}_B(Z(\check U))$ is a subset of 
       ${\cal C}[{\cal A}]^{W_{\Theta}.}$. Set $R={\cal 
       P}_B(Z(\check U))$ and recall that $R$ is  
       generated by the set $\{{\cal P}_B(z_{2\omega_i})|1\leq i\leq n\}.$  
       It is sufficient to show that $\hat m(2\eta_j)\in R$ for 
       all $1\leq j\leq n$. We do this by induction on $j$. It 
       follows from  Lemma 
       3.5 that $\tilde\omega_1=\eta_1$.    Hence Lemma 6.2 and Lemma 
       6.4 ensure that  
       ${\cal P}_B(z_{2\omega_1})$ is a nonzero scalar multiple of $\hat 
       m(2\eta_1)$.    Therefore, $\hat m(2\eta_1)$ is in $R$.   Now assume 
       that $\hat m(2\eta_k)\in R$ for all $1\leq k<j$.   Assume first 
       that $j$ is even and write $j=2l$.   By the inductive 
       hypothesis, 
       $\hat m(2\eta_{l-i})$ and $\hat m(2\eta_{l+i})$ are both in $R$ for 
       $1\leq i<l$.  Hence $R$ contains 
       $\hat m(2\eta_{l-i})\hat m(2\eta_{l+i})$ for $1\leq i<l$.   Now $R$ 
       also contains ${\cal P}_B(z_{2\omega_{2l}})$.    Thus by Lemma 
       6.7(i), $R$ contains 
       $\hat m(2\eta_{2l})$.   The case for j odd is similar using 
       Lemma 6.7(ii) instead of Lemma 6.7(i). $\Box$
       
       \medskip
       
   We now turn our attention  to   case where 
   $\gg,\gg^{\theta}$ is not of type AII.    It 
   follows from Araki's classification [A] and the list in Section 3
  ((3.11)-(3.15)) that 
   $\Sigma$ is a root system of type $B_t$ or $C_t$.  
    Let $\pi_{\gr}$ denote the subset of simple roots in $\pi$ which 
   span the root system of $\gr$. Since the restricted root system of 
   $\gr,\gr^{\theta}$ has rank $t-1$ and $\gr,\gr^{\theta}$  is of type AII, it follows that $\pi_{\gr}$ 
  has    $2t-1$ elements. Moreover, using the standard numbering of 
  the simple roots in Dynkin diagrams found in [H], we have that 
  $\pi_{\gr}=\{\alpha_1,\dots, \alpha_{2t-1}\}$.       Set $W_{\gr}$ equal to the Weyl group of 
  the root system of $\gr$   and note that $W_{\gr}$ is a subgroup of $W$.

  We use below the notation introduced before   Lemma 6.6 where 
    $\Sigma_s'=\{\mu_1,\dots, \mu_{t-1}\}$.  Note that 
  $\Sigma_s\setminus \Sigma_s'=\{\mu_t\}$.  
  Hence $${\cal N}=z^{-\mu_t}\sum_{\gamma\in 
  Q^+(\Sigma)}z^{-\gamma}.\eqno{(6.9)}$$  The next lemma gives technical 
  information which relates certain restricted roots of 
  $\gg,\gg^{\theta}$ to that of $\gr,\gr^{\theta}$.  
  
  \begin{lemma}   Suppose that 
       $\gg,\gg^{\theta}$ is not of type AII.   Then for all $i,k,$ and $s$ such 
   that  $1\leq r\leq k\leq t-1$ and $0\leq s\leq min(r,t-1-k)$
   we have  
  $$\eta_{r-s}+\eta_{k+s}-\eta'_{r-s}-\eta'_{k+s}=\eta_r+\eta_k- 
      \eta_r'-\eta_k'.\eqno{(6.10)}$$
      Moreover, if $\Sigma$ is of type $C_t$, then (6.10) also holds 
     for $s=min(r,t-k)$. 
      \end{lemma}

  \noindent
  {\bf Proof:} Suppose that $r$ has been chosen so that $1\leq r\leq t-1$. We have
 $(\eta_r-\eta_r',\alpha_j)=0$ for all $\alpha_j\in 
 \Sigma_s'$.  Since $\Sigma'$ is of type $A_{t-1}$, it follows  that  the coefficient of $\mu_{t-1}$ in $\eta_r'$ is 
 $r/t$ (see [H, Section 12.2]).  Hence 
 $$(\eta_r-\eta_r',\alpha_t)=-r (\mu_{t-1},\alpha_t)/t.\eqno{(6.11)}$$  Assume first 
 that $\Sigma$ is of type $B_t$. It follows that 
 $\eta_r-\eta_r'=2r\eta_{t}/t$. Hence  $\eta_r+\eta_k-\eta_r'-\eta_k'=2(r+k)\eta_t/t$ for $1\leq 
 r\leq k\leq t-1$ with $0\leq s\leq min(r,t-1-k)$.  In particular, 
 $\eta_r+\eta_k-\eta_r'-\eta_k'$ only depends on $r+k$.   This proves 
 the lemma when $\Sigma$ is of type $B_t$.
 Now assume that  $\Sigma$ is of type $C_t$, 
It follows from (6.11) that   $\eta_r-\eta_r'=r\eta_t/t$ for all 
$1\leq r\leq t-1$.   Now $\eta_t'=0$ so $\eta_r-\eta_r'=r\eta_t/t$ is 
also true for $r=t$.   It follows  that for $1\leq 
 r\leq k\leq t$ with $0\leq s\leq min(r,t-k)$, we have 
 $\eta_r+\eta_k-\eta_r'-\eta_k'=(r+k)\eta_t/t$. Thus 
$\eta_r+\eta_k-\eta_r'-\eta_k'$ depends only on $r+k$ which proves 
the lemma.
 $\Box$

 \medskip
 Recall [Jo, Lemma 7.1.17 and 7.1.25] that the ordinary 
   Harish-Chandra map ${\cal P}$ defines an isomorphism from   $Z(\check 
   U)$ onto
   ${\cal C}[\tau(2\lambda)|\lambda\in P(\pi)]^{W.}$.   Furthermore, 
   the elements 
   $\hat \tau(2\gamma)$ where $\gamma\in P^+(\pi)$ defined by (5.2) are contained in 
   ${\cal C}[\tau(2\lambda)|\lambda\in P(\pi)]^{W.}$. For each $i$ 
   satisfying $1\leq i\leq n$,   let  $z_i$ be the element in 
   $Z(\check U)$ such that 
 ${\cal P}(z_i)=\hat \tau(2\omega_i)$.
 The following is an extended version of Lemma 6.4.
   
   \begin{lemma}  Suppose that 
       $\gg,\gg^{\theta}$ is not of type AII.   Then for each integer 
       $i$ 
   satisfying 
   $2\leq 2i\leq t$, there exist 
   Laurent polynomials $g_{2s}(q)$, $0\leq s\leq i$,  such that  
   $$\eqalign{{\cal P}_B(z_{2i})&=\sum_{0\leq s\leq i}g_{2s}(q)\hat 
   m(2(\eta_{i-s}+\eta_{i+s}))\cr&+\sum_{\gamma\in Q^+(\Sigma)}{\bf 
   Q}(q)\tau(2(2\eta_i-\mu_t-\gamma)) \cr}$$ up to a nonzero scalar.
Similarly, for each integer $i$ such that $1\leq 2i+1\leq t+1$, there 
exist Laurent polynomials $g_{2s+1}(q)$, $0\leq s\leq i$ such that 
    $$\eqalign{{\cal P}_B(z_{2i+1})&=\sum_{0\leq s\leq i}g_{2s+1}(q)\hat 
   m(2(\eta_{i-s}+\eta_{i+1+s}))\cr&+\sum_{\gamma\in Q^+(\Sigma)}{\bf 
   Q}(q)\tau(2(\eta_i+\eta_{i+1}- \mu_t- \gamma)) \cr}$$ up to a 
   nonzero scalar. Moreover, the coefficients of the powers of $q$ in each $g_j(q)$ are 
   nonnegative integers and 
   $g_j(1)=2^j$ for $0\leq j\leq 2i+1$.
   \end{lemma}

   \noindent
   {\bf Proof:}   
   Let $\psi$ denote the isomorphism from ${\cal 
   C}[\tau(2\gamma)|\gamma\in P(\pi)]$ to ${\cal C}[P(\pi)]$ 
   defined by $\psi(\tau(2\gamma)q^{(\rho,2\gamma)})=z^{\gamma}$ 
   for all $\gamma\in P(\pi)$.   Note that for each $i$, 
    $\psi({\cal P}(z_i))=\sum_{\beta\in W\omega_i}z^{\beta} $ up to a 
    nonzero scalar.
   Furthermore,  $$\sum_{\beta\in 
   W\omega_i}z^{\beta}=\sum_{\beta\in 
   W_{\gr}\omega_i}z^{\beta}+\sum_{\gamma\in W\omega_i\setminus 
   W_{\gr}\omega_i} z^{\gamma}.\eqno{(6.12)}$$
  
    Suppose that 
     $\gamma\in W\omega_i\setminus 
   W_{\gr}\omega_i$.   It follows that there exists a reflection 
   $s=s_{\alpha_j}$ 
   corresponding to a simple root $\alpha_j$ in $\pi\setminus \pi_{\gr}$ and 
   elements $w_1\in W_{\gr}$ and $w_2\in W$ such that 
   $\gamma=w_2sw_1\omega_i$.   We may assume that 
   $sw_1\omega_i\neq w_1\omega_i$ and thus  $(\alpha_j,w_1\omega_i)$ is nonzero.
We may further assume that there is a reduced expression 
   for $w_2s_jw_1$ equal to  the product of the reduced expression 
   for $w_2$, $s$, and  the reduced expression for $w_1$ in that order.    
   It follows that  
     $$\omega_i-\gamma-\alpha_j\in Q^+(\pi).\eqno{(6.13)}$$   Now $(\alpha_j,\omega_i)=0$ 
   since $\alpha_i\in \pi_{\gr}$.   Furthermore, 
   $\omega_i-w_1\omega_i \in Q^+(\pi_{\gr})$ because $w_1\in W_{\gr}$.   Thus 
   $(\alpha_j,\alpha)$ must be nonzero for some $\alpha\in 
   \pi_{\gr}$.    We can write $\pi_{\gr}=\{\alpha_1,\dots, 
   \alpha_k\}$ for some choice of $k$ with $1\leq k\leq n$.    
   Checking the   the list (3.11)-(3.15), one sees 
   that $\Delta$ is a root system of type $A_n$, $C_n$, or $D_n$. 
   Hence $j$ 
   must equal $k+1$ and moreover, $\tilde\alpha_j$ is the last simple 
   root of $\Sigma$.   In particular, $k=t-1$ and 
   $\tilde\alpha_j=\mu_t$.   Hence (6.13) implies that $
   \tilde\omega_i-\mu_t-\tilde\gamma\in Q^+(\Sigma)$.   Thus (6.12) 
   yields 
   $$\sum_{\beta\in W\omega_i}z^{\tilde{\beta}}\in \sum_{\beta\in 
   W_{\gr}\omega_i}z^{\tilde{\beta}}+z^{\tilde{\omega_i}}{\cal N}.$$   
   The first assertion of the 
   lemma now follows using (6.4),(6.5), and      Lemma 6.4.
   A similar argument works for the second assertion.
   $\Box$
   
  \medskip
     The next lemma provides additional  
  information about the image of certain central elements under ${\cal 
  P}_B$ not covered by Lemma 6.10.
  
  \begin{lemma}   The element  $\hat 
   m(2\eta_1)$ is in ${\cal P}_B(Z(\check U))$. Furthermore, if
     $\Sigma$ is of type $B_t$ 
 then $\hat 
   m(2\eta_t)$ is in ${\cal P}_B(Z(\check U))$. 
  \end{lemma}
  
  \noindent
  {\bf Proof:}  Fix $i$ with $1\leq i\leq t$.   Suppose that   $\eta_i$ is a minuscule or 
  pseudominuscule fundamental weight in $P(\Sigma)$ (see [M] or [L5, Section 
  7]) and that $\tilde\omega_i=\eta_i$. In particular, the only possible weight less than $\eta_i$ is 
  $0$
  (see [L5, (7.2) and (7.3)]).  By the discussion before Lemma 6.4, we can write ${\cal 
  P}_B(z_{2\omega_i})$ as a linear combination of the $\hat m(2\gamma)$ 
  with $\gamma\in P^+(\Sigma)$ and $\gamma<\eta_i$.
  Hence there exist scalars $a$ and $b$ with $a$ nonzero such that
  ${\cal P}(z_{2\omega_i})=a\hat m(2\eta_i)+b$. It follows that
 $\hat m(2\eta_i)\in {\cal P}_B(Z(\check U))$.  
 
 Since $\Sigma$ is of type A$_t$, B$_t$, or C$_t$, it follows that 
 $\eta_1$ is a minuscule or pseudominuscule weight in $P(\Sigma)$. 
 Moreover, Lemma 3.5(ii) ensures that $\tilde\omega_1=\eta_1$.  Hence $\hat 
 m(2\eta_1)\in {\cal P}_B(Z(\check U))$.  Now assume that $\Sigma$ is of 
 type $B_t$. By Lemma 3.5(ii), we have $\tilde\omega_n=\eta_t$. 
 The fact that
  $\Sigma$ is of type $B_t$ ensures that $\eta_t$ is a minuscule 
 fundamental weight in $P(\Sigma)$. Thus $\hat 
   m(2\eta_t)$ is also in ${\cal P}_B(Z(\check U))$. 
  $\Box$

  \medskip
 For each $1\leq k\leq n$, let $Z_k$ be the image of ${\cal P}_B(z_{k})$ in ${\cal C}[P(\Sigma)]$
  using the isomorphism discussed before Lemma 6.5. In particular, 
  given   integers $l,k$ such that $2\leq 2l\leq t$ and $1\leq 
  2k+1\leq t$, we have
  $$Z_{2l}\in \sum_{0\leq s\leq l}g_{2s}(q)
  m(\eta_{l-s}+\eta_{l+s})+z^{2\eta_l}{\cal N}$$ and $$Z_{2l+1}\in \sum_{0\leq s\leq 
  l}g_{2s+1}(q) 
  m(\eta_{k-s}+\eta_{k+1+s})+z^{\eta_k+\eta_{k+1}}{\cal N} $$ where the 
  $g_j$ are   Laurent polynomials in $q$ satisfying $g_j(1)=2^j$. 
  Set $$Z_{2l}'=\sum_{0\leq s\leq l}g_{2s}(q)
  m'(\eta'_{l-s}+\eta'_{l+s})$$ and $$Z'_{2l+1}= \sum_{0\leq s\leq 
  l}g_{2s+1}(q) 
  m'(\eta'_{k-s}+\eta'_{k+1+s}). $$

  Assume for the moment that  ${\cal C}[{\cal A}]^{W_{\Theta}.}$ is a subset of 
  ${\cal P}_B(\check U^B)$.  Lemma 4.3 and (5.5) imply that 
  $$\eqalign{\{{\rm tip}(u)|\ u\in {\cal P}_B(\check 
  U^B)\}=&\{{\rm tip}(u)|\ u\in {\cal P}_B(Z(\check U))\}\cr=&{\rm 
  span}\{\tau(-2\tilde\mu)|\ \mu\in P^+(\Sigma)\}.\cr}$$  Arguing as 
  in the proof of Theorem 5.4, it then follows that 
  ${\cal P}_B(\check U^B)={\cal P}_B(Z(\check U))={\cal C}[{\cal A}]^{W_{\Theta}.}$. 
 Hence  the  following generalization of Theorem 6.8 completes the proof of 
Theorem 6.1.

    \begin{theorem} Suppose $\gg,\gg^{\theta}$ contains a symmetric pair $\gr,\gr^{\theta}$ 
   of type AII
   and the rank of the restricted root system associated to 
   $\gr,\gr^{\theta}$ is equal to $t-1$ with $t-1\geq 3$.  Then
     the  ${\cal C}[{\cal A}]^{W_{\Theta}.}$ is a subalgebra of 
     ${\cal P}_B(Z(\check U))$.
       \end{theorem}
    
  \noindent
  {\bf Proof:} Let  $R$ be the image of ${\cal P}_B(Z(\check U))$ 
  using the isomorphism described before Lemma 6.5.  In particular, 
  $R$ contains $Z_k$ for $1\leq k\leq n$.
 It is sufficient to show that $ m(\eta_i)\in R$ for 
  $1\leq i\leq t$.  By Theorem 6.8, the result is true when $\Sigma$ is of type 
  $A_t$.   Hence we may assume that $\Sigma$ is of type $B_t$ or $C_t$.
  By Lemma 6.11, $R$ contains $m(\eta_1)$.   
 We use
  induction   along 
  the lines of  the proof of Theorem 6.8 to show that all the $m(\eta_{j})$ 
 are  in $R$.
  
  Set $t'=t-1$ if $\Sigma$ is of type $B_t$ and $t'=t$ if $\Sigma$ is 
  of type $C_t$.  For each $1\leq j\leq t-1$, 
 we have that the restriction $\eta'_j$ of $\eta_j$ to the root 
 system  $\Sigma'$   is just
  the fundamental weight in $P^+(\Sigma')$ associated to the root 
  $\mu_j$.   By (6.4) and Lemma 6.9, we have that 
   $$m(\eta_{i-s}+\eta_{k+s})\in 
  z^{\eta_i+\eta_k-\eta_i'-\eta_k'}m'(\eta'_{i-s}+\eta_{k+s}')
 +z^{\eta_i+\eta_k}{\cal N}  $$ for all $i, k, $ and $s$ such that 
 $1\leq i\leq k\leq t-1$, $i+k\leq t'$, and $0\leq s\leq i$.
 It follows from the definition of the $Z_k$ and $Z_k'$    that 
 $$Z_{2l}\in 
  z^{2\eta_l-2\eta_l'}Z'_{2l}
 +z^{2\eta_l}{\cal N}$$ for all integers $l$ such that $2\leq 2l\leq t$
 and $$Z_{2j+1}\in 
  z^{ \eta_j+\eta_{j+1}- \eta_j'-\eta_{j+1}'}Z'_{2j+1}
 +z^{\eta_j+\eta_{j+1}}{\cal N}$$
for
 all integers $j$ such that $1\leq 2j+1\leq t$. 
By (6.5) and Lemma 6.9, we further have that $$m(\eta_{i-s})m(\eta_{k+s})\in 
  z^{\eta_i+\eta_k-\eta_i'-\eta_k'}m'(\eta'_{i-s})m'(\eta_{k+s}')
 +z^{\eta_i+\eta_k}{\cal N} $$  for all $i,k,$ and $s$ such that 
 $1\leq i,k\leq t-1$, $i+k\leq t'$, and $0\leq s\leq i$.

   Assume that $m(\eta_1),\dots, m(\eta_{s})$ are in $R$ for some $1\leq 
  s\leq t'-1$.   Suppose that $$s+1=2l{\rm \ for \ some \ integer \ 
  }l.\eqno{(6.14)}$$   
Arguing as in Lemma 6.7, we can write 
  $m'(\eta'_{2l})$ as a linear combination of $Z'_{2l}$ and 
  the elements in the set $\{m'(\eta'_{l-j})m'(\eta'_{l+j})|0\leq 
  j<l\}$.  Thus there is an element $Y$ in $R$ 
  which is a  linear 
  combination of elements in the set $\{Z_l\}\cup\{m(\eta_{l-j})m(\eta_{l+j})|0\leq 
  j<l\}$ such that $$Y\in m(\eta_{2l})+z^{2\eta_{l}}{\cal 
  N}.$$ 
 Consider an element  $\gamma\in P^+(\Sigma)$ such that
  $z^{\gamma}\in z^{2\eta_{l}}{\cal 
  N}$.  By (6.9), we have that $2\eta_l-\mu_t-\gamma\in Q^+(\Sigma)$.  
  
  Suppose that $\Sigma$ is of type $B_t$. By Lemma 6.11, both 
  $m(\eta_1)$ and $m(\eta_t)$ are in $R$. It follows from the list of the fundamental weights in [H, 
  Section 13.2]  that the coefficient of $\mu_t$ in $2\eta_l$ is $2l$. 
  Hence the 
  coefficient of $\mu_t$ in $\gamma$ must be a nonnegative integer
  less than or equal to $2l-1$. Recall (6.14) that $2l-1=s$.   The table  in ([H, Section 
  13.2]) yields that $\gamma$ is a 
  linear combination of elements in the set $\{\eta_i|\ 1\leq 
  i\leq s$ or $i=t\}$.   Hence $Y-m(\eta_{2l})$ is in the subalgebra of $R$ 
  generated by the set $\{m(\eta_i)|\   1\leq i\leq s$ or $i=t\}$.   
  The fact that $Y\in R$ forces $m(\eta_{2l})\in R$.   A similar argument shows that 
  $m(\eta_{2l+1})$ is  in $R$ when $s+1=2l+1$ for some integer $l$.  
  Recall that we have assumed that $s\leq t'-1=t-2$.  Thus by   induction, $R$ contains 
  $m(\eta_2),\dots, m(\eta_{t-1})$ which completes 
  the proof of the theorem for $\Sigma$  of type $B_t$.

 Now assume that $\Sigma$ is of type $C_t$. By Lemma 6.11, $R$ 
 contains $m(\eta_1)$. Assume that $s+1=2l$ for some integer $l$. It follows from [H, 
 Section 13.2, Table 1] that  the coefficient of $\mu_t$ in $2\eta_l$ is $l$. 
Hence the
  coefficient of $\mu_t$ in $\gamma$ must be a nonnegative integer
  less than or equal to $l-1$.  Note that $2l=s+1\leq t$. Hence 
  $l-1\leq (t-2)/2$. Checking the table in [H, Section 
  13.2], we see that $\gamma$ is a linear combination of weights in 
  the set $\{\eta_i|\ 1\leq 
  i\leq t-2\}$.  Hence $Y-m(\eta_{2l})$ is in the subalgebra of $R$ 
  generated by the set $\{m(\eta_i)|\   1\leq i\leq t-2\}$.  
    Since $Y\in R$, we must 
  also have $m(\eta_{2l})\in R$.  A similar argument shows that
  $m(\eta_{l}+\eta_{l+1})$ is in $R$ for  $s=2l+1$.     Thus by induction, $R$ contains 
  $m(\eta_1),\dots, m(\eta_{t})$.  $\Box$

  \section{Appendix: Commonly used  notation}

\noindent
Here is a list of notation defined in Section 1 (in the following 
order):

\noindent
${\bf C}$, ${\bf Q}$,
${\bf Z}$, ${\bf R}$,
${\bf N}$, $q$,
${\cal C}$, ${\cal R}$,
$Q(\Phi)$,
$Q^+(\Phi)$,
$P(\Phi)$
$P^+(\Phi)$,
$\gg$,
$\gn^-\oplus \gh \oplus \gn^+ $,
$\theta$,
$\gg^{\theta}$
$\Delta$,
$(\ ,\ )$,
$\pi=\{\alpha_1,\dots, \alpha_n\}$,
$\leq$,
$\Sigma$,
$\Theta$,
$\tilde\alpha$,
${\rm p}(i)$,
$\pi^*$,
$U$, $U_q(\gg)$,
$x_i$, $y_i$, $t_i^{\pm 1}$,
$T$,
$U^+$,
$G^-$, $U^0$, $\tau$, $S_{\beta}$,
$A_+$,
$\check U$,
$\check T$, $\check U^0$,
$\pi_{\Theta}$
${\cal M}$,
$T_{\Theta}$,
$\tilde\theta$,
$B_i$,
$B_{\theta}$,
${\cal S}$, ${\cal D}$, $U_{\cal R}$,
$ {\cal B}$,
${\bf H}$,
$F_r(\check U)$,
$\adr$,
$\check U^B$, $Z(\check U)$, $z^{\lambda}$, 
$L(\Lambda)$,
${\cal C}[G]$ ($G$ a multiplicative group), ${\cal C}[H]$ ($H$ an 
additive group).
\medskip
\begin{tabbing}
	\noindent
\=Defined in Section 2:\=\\
\>${\cal M}^+$\> ${\cal M}\cap U^+$\\
\>${\cal M}^-$\> ${\cal M}\cap G^-$\\
\>$\ad $\> left adjoint action\\
\>$N^+$\>  subalgebra  of $U^+$  generated  by $(\ad {\cal M}^+){\cal 
C}[x_i|\ \alpha_i\notin\pi_{\Theta}\}$\\
\>$N^-$\>subalgebra  of $G^-$  generated  by $(\ad{\cal M}^-){\cal 
C}[y_it_i|\alpha_i\notin\pi_{\Theta}]$\\
\>$S_{\beta,r}$\>the  sum  of  weight spaces $S_{\beta'}$  with
$\tilde\beta'=\tilde\beta$\\
\>$T'_{\geq}$\>monoid  generated by $t_i^2$, $\alpha_i\in \pi^*$, 
and $t_i$, $\alpha_i\in {\cal S}$\\
\>$T'$\>group  generated by $t_i$ for $\alpha_i\in \pi^*$\\
\>$\check {\cal A}$\> $\{\tau(\tilde\mu)|\ \mu\in P(\pi)\}$\\
\>$\check T_{\Theta}$\> $\{\tau(1/2(\mu+\Theta(\mu))|\ \mu\in 
P(\pi)\}$\\
\>${\cal P}_B$\> see Definition 2.3\\
\>${\cal A}$\>$\{\tau(2\mu)|\ \mu\in P(\Sigma)\}$\\
\>$G^+$\>  algebra   generated by $x_it_i^{-1}$ for $1\leq 
i\leq n$\\
\>$U^-$\> algebra generated by $y_i$ for $1\leq i\leq n$\\
\>\>\\
\=Defined in Section 3:\=\\
\>$\omega_i$\>   fundamental weight corresponding to $\alpha_i$\\
\>$\omega_i'$\>  fundamental restricted weight corresponding to 
$\tilde\alpha_i$\\
\>$\pi_i$\>$\{\alpha_j|\ (\omega_j,\tilde\alpha_i)\neq 0\}$\\
\>$\gg_i$\> semisimple Lie algebra with root system $\pi_i$\\
\>$\Sigma_i$\> restricted roots of $\gg_i,\gg_i^{\theta}$\\
\>$W$\>   Weyl group of the root system of $\gg$\\
\>$W_{\Theta}$\>   Weyl group of the root system $\Sigma$\\
\>\>\\
\=Defined in Section 4:\=\\
\>${\cal X}$\> the radial component map (see Theorem 4.1)\\
\>${\cal C}[Q(\Sigma)]{\cal A}$\> subring of ${\rm End}_r{\cal 
C}[P(\Sigma)]$ generated by ${\cal C}[Q(\Sigma)]$ and ${\cal A}$\\
\>${\cal C}(Q(\Sigma)){\cal A}$\> localization of ${\cal 
C}[Q(\Sigma)]{\cal A}$ at ${\cal C}[Q(\Sigma)]\setminus\{0\}$\\
\>$g_{\lambda}$\> zonal spherical function at $\lambda$\\
\>$\varphi_{\lambda}$\> image of $g_{\lambda}$ in ${\cal 
C}[P(2\Sigma)]$\\
\>${\cal F}$\> filtration on $\check U$ defined by (4.1), (4.2), and 
(4.3)\\
\>$p$\> see (4.4)\\
\>tip($X$)\> highest degree homogeneous term of $X$\\
\>\>\\
\=Defined in Section 5:\=\\
\>${\cal P}$\> the  Harish-Chandra  map 
defined using (5.1)\\
\>$z_{2\mu}$\>  unique central element in $\tau(2\mu)+(\adr 
U_+)\tau(2\mu)$\\
\>$\rho$\> half sum  of the positive roots in $\Delta$\\
\>$\hat\tau(\lambda)$\> $\sum_{w\in 
W}\tau(w\lambda)q^{(\rho,w\lambda)}$\\
\>$a_{2\mu}$\> scalar defined in Lemma 5.1\\
\>$w_0$\> longest element of $W$\\
\>$w.q^{(\rho,mu)}\tau(\mu)$\> $q^{(\rho,w\mu)}\tau(w\mu)$\\
\>$W'$\>   Weyl group associated to root system of $\pi_{\Theta}$\\
\>$\hat m(2\tilde\mu)$\>$\sum_{\gamma\in 
W_{\Theta}\tilde\mu}q^{(\rho,2\gamma)}\tau(2\gamma)$\\
\>$\check U_{{\bf C}(q)}$\> ${\bf C}(q)$ algebra generated by 
$x_i,y_i,1\leq i\leq n,$ and $\check T$\\
\>$A$\> ${\bf C}[q]_{(q-1)}$\\
\>$\hat U$\> $A$ algebra generated by $x_i,y_i$, $1\leq i\leq n$, 
 ${{(t-1)}\over{(q-1)}}$, $t\in \check T$\\
 \>$h_1,\dots, h_n$\> basis for $\gh$\\
 \>$T_2$\>$\{\tau(2\mu)|\ \mu\in P(\pi)\}$\\
\>\>\\
\=Defined in Section 6:\=\\
\>$t$\> The rank of $\Sigma$\\
\>$\gr,\gr^{\theta}$\> proper maximal symmetric pair in 
$\gg,\gg^{\theta}$ of type AII\\
\>$\mu_1,\dots, \mu_t$\>    simple positive roots of $\Sigma$\\
\>$\eta_1,\dots,\eta_t$\>    restricted fundamental weights\\
\>$\eta_0$\>$0$\\
\>$\eta_{t+1}$\>$0$\\
\>$s_{\alpha}$\>  reflection in $W$ corresponding to root $\alpha$  \\
\>$m(\lambda)$\>$\sum_{\gamma\in W_{\Theta}\lambda}z^{\gamma}$\\
\>$\Sigma_s'$\> a subset of $\{\mu_1,\dots, \mu_t\}$\\
\>$\Sigma'$\> root system generated by $\Sigma_s'$\\
\>$\lambda'$\> $\lambda'\in P(\Sigma')$ and
$(\lambda-\lambda',Q(\Sigma'))=0$\\
\>$W_{\Theta}'$\> Weyl group of $\Sigma'$\\
\>${\cal N}$\>$\sum_{\mu_i\in \Sigma_s\setminus 
\Sigma'_s}\sum_{\gamma\in Q^+(\Sigma)}{\bf N}z^{-\mu_i-\gamma}$\\
\>$\hat M_{k}$\> see (6.7) and (6.8)\\
\>$M_k$\> the image of $\hat M_k$ in ${\cal C}[P(\Sigma)]$\\
\>$\pi_{\gr}$\> the  set of simple roots of the root system of $\gr$\\
\>$W_{\gr}$\> the Weyl group of the root system of $\gr$\\
\>$z_i$\> $z_i\in Z(\check U)$ and ${\cal P}(z_i)=\hat \tau(2\omega_i)$\\
\>$Z_k$\> the image of ${\cal P}_B(z_k)$ in ${\cal C}[P(\Sigma)]$\\
\>$Z_k'$\> see the discussion before Theorem 6.12\\
\end{tabbing}

\medskip

 \centerline{REFERENCES}

\bigskip
\bigskip
\noindent   [A] S. Araki, On root systems and an infinitesimal 
classification of irreducible symmetric spaces, {\it Journal of 
Mathematics, Osaka City University} {\bf 13} (1962), no. 1, 1-34.

\medskip
\noindent [DN] M.S. Dijkhuizen and M. Noumi, A family of quantum
projective spaces and related $q$-hypergeometric orthogonal
polynomials, {\it Transactions of the A.M.S.} {\bf 350} (1998), no.
8, 3269-3296.

\medskip
\noindent [DS] M.S. Dikhhuizen and J.V. Stokman, Some limit transitions between BC type orthogonal polynomials interpreted on
quantum complex Grassmannians, {\it Publ. Res. Inst. Math. Sci. 35} (1999), 451-500.

\medskip\noindent [H]   J.E. Humphreys, {\it Introduction to Lie Algebras
and Representation Theory}, Springer-Verlag, New York
(1972).

\medskip\noindent [He] S. Helgason, Some results on invariant 
differential operators on symmetric spaces, {\it American Journal of 
Mathematics} {\bf 114}, No. 4, 789-811. 

\medskip\noindent[HC1] Harish-Chandra, On some applications of the 
universal enveloping algebra of a semisimple Lie group, {\it 
Transactions of the American Mathematical Society} {\bf 70} (1951) 
28-96.

\medskip\noindent[HC2] Harish-Chandra, Spherical functions on a 
semisimple Lie group, I, {\it American Journal of Mathematics}  {\bf 80} 
(1958) 241-310.

\medskip\noindent [Ja] N. Jacobson, {\it Basic Algebra I,} W. H. 
Freeman and Company, San Francisco (1974).
 
\medskip \noindent [JL] A. Joseph and G. Letzter, Local
finiteness of the adjoint action for quantized enveloping
algebras, {\it Journal of Algebra} {\bf 153} (1992), 289 -318.

\medskip
\noindent [JL2] A. Joseph and G. Letzter, Separation of variables
for quantized enveloping algebras, {\it American Journal of
Mathematics} {\bf 116} (1994), 127-177.
 
\medskip\noindent  [Jo] A. Joseph, {\it Quantum Groups and Their Primitive
Ideals}, Springer-Verlag, New York (1995).

\medskip
\noindent [Ke] M.S. K\'eb\'e, ${\cal O}$-alg\`ebres quantiques, {\it C. R. Acad\'emie des Sciences,
S\'erie I, Math\'ematique } {\bf 322} (1996), no. 1, 1-4.

\medskip
\noindent
[Kn] A. W. Knapp, {\it Lie Groups Beyond an Introduction}, Progress 
in Math. {\bf 140}, Birkh\"auser, Boston (1996).

\medskip\noindent [Le] J. Lepowsky, On the Harish-Chandra Homomorphism,
{\it Transactions of the American Mathematical Society} {\bf 208} 
(1975) 193-218.

\medskip\noindent [L1] G. Letzter, Symmetric pairs for quantized 
enveloping algebras, {\it Journal of Algebra} {\bf 220} (1999), no. 2, 
729-767.

\medskip \noindent [L2] G. Letzter, Harish-Chandra modules for
quantum symmetric pairs, {\it Representation Theory, An Electronic
Journal of the AMS} {\bf 4} (1999)
64-96.

\medskip \noindent [L3] G. Letzter, Coideal subalgebras and quantum 
symmetric pairs, In: {\it New Directions in Hopf Algebras, MSRI publications}
{\bf 43}, Cambridge University Press (2002), 117-166.

\medskip\noindent [L4]  G. Letzter, Quantum symmetric pairs and their 
zonal spherical functions, {\it Transformation Groups} {\bf 8} 
(2003), no. 3. 261-292.

\medskip\noindent [L5]  G. Letzter, Quantum  zonal spherical  
functions and Macdonald polynomials, Advances in Mathematics, in
press (corrected proof available online).

\medskip\noindent [L6]  G. Letzter, Invariant differential operators 
for quantum symmetric spaces, II, preprint (math.QA/0406194).

\medskip
\noindent
[M] I.G. Macdonald,
Orthogonal polynomials associated with root systems, {\it S\'eminaire Lotharingien de Combinatoire} 
{\bf 45} (2000/01) 40 pp.

\medskip
\noindent [N] M. Noumi, Macdonald's symmetric polynomials as zonal
spherical functions on some quantum homogeneous spaces, {\it
Advances in Mathematics} {\bf 123} (1996), no. 1, 16-77.

\medskip
\noindent [NDS] M. Noumi, M.S. Dijkhuizen, and T. Sugitani, Multivariable Askey-Wilson polynomials and quantum complex
Grassmannians, {\it Fields Institute Communications} {\bf 14} (1997), 167-177.

\medskip
\noindent [NS] M. Noumi and T. Sugitani, Quantum symmetric spaces
and related q-orthogonal polynomials, in: {\it Group Theoretical
Methods in Physics (ICGTMP)} (Toyonaka, Japan, 1994), World
Science Publishing, River Edge, New Jersey (1995), 28-40.

\end{document}